\documentclass[reqno, a4]{amsart}

\usepackage{harvard}
\usepackage{epsfig}                     
\usepackage{ifthen}

\usepackage{amsmath,amssymb,amsthm}
\usepackage{amsbsy,bm}                  
\usepackage[english]{babel}
\usepackage[latin1]{inputenc}           
\usepackage{fontenc}
\usepackage{url}
\usepackage{bbm}
\usepackage[colorlinks=true]{hyperref}  

\let\margin\marginpar
\newcommand\myMargin[1]{\margin{\raggedright\scriptsize #1}}
\renewcommand{\marginpar}[1]{\myMargin{\textsf{#1}}}

\renewcommand{\marginpar}[1]{}

\theoremstyle{plain}                    
\begingroup                             
\newtheorem{thm}{Theorem}[section]

\newtheorem{cor}[thm]{Corollary}
\newtheorem{lem}[thm]{Lemma}
\newtheorem{prop}[thm]{Proposition}

\endgroup

\theoremstyle{definition}
\newtheorem{defn}[thm]{Definition}

\theoremstyle{remark}
\newtheorem{rem}[thm]{Remark}           

\numberwithin{equation}{section}

\newcommand{\nothing}[1]{}
\newcommand{\bs}[1]{\boldsymbol{#1}}
\def\half{\frac{1}{2}}

\newcommand{\eps}{\varepsilon}

\DeclareMathOperator{\pr}{pr}
\DeclareMathOperator{\vc}{vc}
\DeclareMathOperator{\sign}{sign}
\DeclareMathOperator{\card}{card}

\def\argmin{\operatornamewithlimits{arg\,min}}

\begin{document}

\title[Monte Carlo Algorithms for Optimal Stopping]{Monte Carlo Algorithms for Optimal Stopping and Statistical Learning}

\author{Daniel Egloff}
\address{Wasserfuristrasse 42, 8542 Wiesendangen, Switzerland}
\curraddr{Zurich Cantonal Bank, P.O. Box, CH-8010 Zurich, Switzerland. Phone: +41 1 292 45 33}
\email{daniel.egloff@zkb.ch}

\keywords{Optimal stopping, American options,
statistical learning, empirical processes, uniform law of large
numbers, concentration inequalities, Vapnik-Chervonenkis classes,
Monte Carlo methods.}

\date{This Version: \today. First Version: 4th April 2003}

\subjclass[2000]{Primary 91B28, 60G40, 93E20; Secondary 65C05, 93E24, 62G05}

\begin{abstract}
We extend the Longstaff-Schwartz algorithm for approximately
solving optimal stopping problems on high-dimensional state
spaces. We reformulate the optimal stopping problem for Markov
processes in discrete time as a generalized statistical learning
problem. Within this setup we apply deviation inequalities for
suprema of empirical processes to derive consistency criteria, and
to estimate the convergence rate and sample complexity. Our
results strengthen and extend earlier results obtained by
\citeasnoun{Protter-Lamberton:2001}.
\end{abstract}

\maketitle

\section{Introduction}

The problem of arbitrage-free pricing American options has renewed
the interest in efficient methods for numerically solving
high-dimensional optimal stopping problems. In this paper we
explain how to solve a discrete-time, finite-horizon optimal
stopping problem by restating it as a generalized statistical
learning problem. We give a unified treatment of the
Longstaff-Schwartz and the Tsitsiklis-Van Roy algorithm. They use
both Monte Carlo simulation and linearly parameterized
approximation spaces. We introduce a new class of algorithms which
interpolate between the Longstaff-Schwartz and Tsitsiklis-Van Roy
algorithm and relax the linearity assumption of the approximation
spaces.

Learning an optimal stopping rule differs from the standard setup
in statistical and machine learning in the sense that it requires
a series of learning tasks, one for every time step, starting at
the terminal horizon and proceeding backward. The individual
learning tasks are connected by the dynamic programming principle.
At each time step, the result depends on the outcome of the
previous learning tasks. Connecting the subsequent learning tasks
to a recursive sequence of learning problems leads to an error
propagation. We control the error propagation by using a Lipschitz
property and a suitable error decomposition which relies on the
convexity of the approximation spaces. Finally, we estimate the
sample error with exponential tail bounds for the supremum of
empirical processes. To apply these techniques we need to
calculate the covering numbers of certain function classes. An
important type of function classes for which good estimates on the
covering numbers exist are the so called Vapnik-Chervonenkis (VC)
classes, see \citeasnoun{Vaart-Wellner:1996} or
\citeasnoun{Anthony-Bartlett:1999}. We prove that payoff functions
evaluated at Markov stopping times parameterized by a VC-class of
functions is again a VC-class. The covering number estimate of
\citeasnoun{Haussler:1995} then gives the required bounds. Our
approach is conceptually different from
\citeasnoun{Protter-Lamberton:2001}, which is purely tailored to
the classical Longstaff-Schwartz algorithm with linear
approximation. By exploiting convexity and fundamental properties
of VC-classes we can prove convergence and derive error estimates
under less restrictive conditions, also if both the dimension of
the approximating spaces and the number of samples tends to
infinity.

This paper is structured as follows. The next background section
discusses recent developments in numerical techniques for optimal
stopping problems and summarizes the probabilistic tools which we
use in this work. Section \ref{sec:intro opt stopping} reviews
discrete-time optimal stopping problems. Section \ref{sec:opt stop
stat learning} shows how to restate optimal stopping as a
statistical learning problem and introduces the dynamic look-ahead
algorithm. In Section \ref{sec:main results} we state and comment
our main results: a general consistency result for convergence,
estimates of the overall error, the convergence rate, and the
sample complexity. The focus of the work lies in estimating the
sample error. The proofs are deferred to Section \ref{sec:proofs}
where we also introduce the necessary tools of the
Vapnik-Chervonenkis theory.

\section{Background}
\label{sec:background}

Optimal stopping problems naturally arise in the context of games
where a player wants to determine when to stop playing a sequence
of games to maximize his expected fortune. The first systematic
theory of optimal stopping emerged with
\citeasnoun{Wald-Wolfowitz:1948} on the sequential probability
ratio test. The monographs by
\citeasnoun{Chow-Robbins-Siegmund:1971} and
\citeasnoun{Shiryayev:1978} provide an extensive treatment of
optimal stopping theory.

The general no-arbitrage valuation of American options in terms of
an optimal stopping problem begins with
\citeasnoun{Bensoussan:1984} and \citeasnoun{Karatzas:1988}.
Nowadays, American option valuation is an important application of
optimal stopping theory. For more background on American options
and financial aspects of the related optimal stopping problem we
refer to \citeasnoun{Karatzas-Shreve-MathFin:1998}.

\subsection{Algorithms for Solving Optimal Stopping Problems}

Optimal stopping problems generally cannot be solved in closed
form. Therefore, several numerical techniques have been developed.
\citeasnoun{Barone-Adesi-Whaley:1987} propose a semi-analytical
approximation. The binomial tree algorithm of
\citeasnoun{Cox-Ross-Rubinstein:1979} directly implements the
dynamic programming principle. Other approaches comprise Markov
chain approximations, see \citeasnoun{Kushner:1998}, direct
integral equation and PDE methods. The PDE methods are based on
variational inequalities, developed in
\citeasnoun{Bensoussan-Lion:1982} or
\citeasnoun{Jaillet-Lamberton-Lapeyre:1990}, the linear
complementary problem, see \citeasnoun{Huang-Pang:1998}, or the
free boundary value problem, see \citeasnoun{VanMoerbeke:1976}.
However, the viability of any of these methods is prohibited by
the curse of dimensionality. For these algorithms the computing
cost and storage needs grow exponentially with the dimension of
the underlying state space.

To address this limitation, new Monte Carlo algorithms have been
proposed. The first landmark papers in this direction are
\citeasnoun{Boessarts:1989}, \citeasnoun{Tilley:1993}, and
\citeasnoun{Broadie-Glass:1997}.
\citeasnoun{Longstaff-Schwartz:1999} introduce a new algorithm for
Bermudan options in discrete time. It combines Monte Carlo
simulation with multivariate function approximation. They show how
to solve the optimal stopping problem algorithmically by a nested
sequence of least-square regression problems and briefly outline a
convergence proof. \citeasnoun{Tsitsiklis-Roy:1999} independently
propose an alternative parametric approximation algorithm on the
basis of temporal--difference learning. Their approach relies on
stochastic approximation of fixed points of contraction maps. They
prove almost sure convergence by using stochastic approximation
techniques as developed in \citeasnoun{Kushner:1978},
\citeasnoun{Benveniste:1990}, or \citeasnoun{Kushner-Yin:1997}.
The Longstaff-Schwartz as well as the Tsitsiklis-Van Roy algorithm
approximate the value function or the early exercise rule and
therefore provide a lower bound for the true optimal stopping
value. \citeasnoun{Rogers:2001} proposes a method based on the
dual problem which results in upper bounds. The overview paper
\citeasnoun{Broadie-Glass:1998} describes the state of development
of Monte Carlo algorithms for optimal stopping as of 1998. A more
recent reference is the book of \citeasnoun{Glasserman:2004}. A
comparative study of various Monte Carlo algorithms for optimal
stopping can be found in \citeasnoun{Madan-Laprise-Wu-Fu:2001}.

Despite of the contributions of \citeasnoun{Tsitsiklis-Roy:1999},
\citeasnoun{Longstaff-Schwartz:1999}, and
\citeasnoun{Rogers:2001}, many aspects of Monte Carlo algorithms
for optimal stopping such as convergence and error estimates
remain unanswered. \citeasnoun{Protter-Lamberton:2001} provide a
complete convergence proof and a Central Limit Theorem for the
Longstaff-Schwartz algorithm. But there are so far no results on
more general possibly nonlinear approximation schemes, the rate of
convergence or error estimates. These problems are the main topics
addressed in this paper.

\subsection{Probabilistic Tools}

The main probabilistic tools which we apply in this paper are
exponential deviation inequalities for suprema of empirical
processes. These tail bounds have been developed by
\citeasnoun{Vapnik-Chervonenkis:1971}, \citeasnoun{Pollard:1984},
\citeasnoun{Talagrand:1994}, \citeasnoun{Ledoux:1996},
\citeasnoun{Massart:2000}, and \citeasnoun{Rio:2001} and many
others. Compared to Central Limit Theorems, they are
non-asymptotic and provide meaningful results already for a finite
sample size. Deviation inequalities together with combinatorial
estimates of covering numbers in terms of the Vapnik-Chervonenkis
dimension are the cornerstones of statistical learning by
empirical risk minimization. For additional details on statistical
learning theory we refer to \citeasnoun{Vapnik:1982},
\citeasnoun{Vidyasagar:1996}, \citeasnoun{Anthony-Bartlett:1999},
\citeasnoun{Vapnik:1999}, \citeasnoun{Cucker-Smale:2001},
\citeasnoun{Mendelson1:2003}, \citeasnoun{Mendelson2:2003}, and
\citeasnoun{Kohler-etal:2002}.

\subsection{Basic Notations}

The following terminology and notation will be used throughout
this paper. If $\mu$ is a measure on a measurable space $(M,
\mathcal{A})$ we denote by $L_p(M,\mu)$ the usual $L_p$-spaces
endowed with the norm $\| \ \|_{p,\mu}$. If we need to indicate
the measure space we write $\| \ \|_{p,M,\mu}$. Let $d_{p,\mu}$ be
the induced metric $d_{p,\mu}(f,g) = \| f-g \|_{p,\mu}$.

Let $(M, d)$ be a metric space. If $U \subset M$ is an arbitrary
subset we define the covering number
\begin{equation}\label{equ:def cov num}
    N(\eps, U, d) = \inf \{n \in \mathbb{N} \mid
        \exists \; \{x_1, \ldots x_n\} \subset M \text{ such that } \forall x\in U \ {\min}_{i=1,\ldots,n} \ d(x,x_i) \leq \eps \}
\end{equation}
which is the minimum number of closed balls of radius $\eps$
required to cover $U$. The logarithm of the covering number is
called the entropy. The growth rate of the entropy for $\eps
\rightarrow 0$ is a measure for the compactness of the metric
space $U$.

Let $X, X_1, X_2, \ldots$ be i.i.d. random elements on a
measurable space $(M, \mathcal{A})$ with distribution $P$. The
empirical measure of a random sample $X_1, \ldots, X_n$ is the
discrete random measure given by
\begin{equation}\label{equ:def empirical measure}
    P_n(A) = \frac{1}{n} \sum_{i=1}^n 1_{\{X_i \in A\}}, \quad A \in
    \mathcal{A},
\end{equation}
or if $g$ is a function on $M$
\begin{equation}\label{equ:def empirical measure func}
    P_n g = \frac{1}{n} \sum_{i=1}^n g(X_i).
\end{equation}
The empirical measure is a random measure supported on
$(M^{\infty}, P^{\infty}, \mathcal{A}^{\infty})$ where $M^{\infty}
= \prod_{\mathbb{N}} M$ is the product space of countably many
copies of $M$, $P^{\infty}$ the product measure, and
$\mathcal{A}^{\infty}$ the product $\sigma$-algebra. The random
variables $X_i$ can now be identified with the $i$-th coordinate
projections.

\section{Review of Discrete Time Optimal Stopping}
\label{sec:intro opt stopping}

Let $\mathbf{X} = (X_t)_{t=0,\ldots,T}$ be a discrete time
$\mathbb{R}^m$-valued Markov process. We assume $\mathbf{X}$ is
canonically defined on the path space $\bs{\mathcal{X}} =
\mathbb{R}^m \times \ldots \times \mathbb{R}^m$ of $T+1$ factors
and identify $X_t$ with the projection onto the factor $t$. We
endow $\bs{\mathcal{X}}$ with the Borel $\sigma$-algebra
$\mathcal{B}$. Let $\mathcal{F}_t$ be the smallest
$\sigma$-algebra generated by $\{X_s \mid s \leq t\}$ and
$\mathbb{F} = (\mathcal{F}_t)_{t = 0,\ldots,T}$ the corresponding
filtration.

Let $P$ be the law of $\mathbf{X}$ on $\bs{\mathcal{X}}$ and
$\mu_t = P_{X_t}$ the law of $X_t$ on $\mathbb{R}^m$. We introduce
the spaces of Markov $L_p$-functions
\begin{equation}\label{equ:space LpX}
    L_p(\mathbf{X}) = \{h = (h_0,\ldots,h_T) \mid h_t \in L_p(\mathbb{R}^m, \mu_t), \; \forall t=0,\ldots,T\},
\end{equation}
with norm
\begin{equation}\label{equ:norm LpX}
    \|h\|_p =
        \sum_{t=0}^T \|h_t\|_{p,\mu_t} = \sum_{t=0}^T E[
        |h_t(X_t)|^p]^{1/p}.
\end{equation}
For brevity we drop the measures $P$, $\mu_t$ and the coordinate
projections $X_t$ in our notation whenever no confusion is
possible. Also, if $h \in L_p(\mathbf{X})$ and
$\mathbf{x}=(x_0,\ldots,x_T) \in \bs{\mathcal{X}}$ is a point of
the path space, we introduce the shorthand notation
\begin{equation}\label{equ:h(xx)}
    h(\mathbf{x})_t \equiv h_t(x_t).
\end{equation}

\subsection{Discrete Time Optimal Stopping}

In the following $f \in L_1(\mathbf{X})$ is a nonnegative reward
or payoff function. The optimal stopping problem consists of
finding the value process
\begin{equation}\label{equ:opt stop problem}
    V_t = {ess \sup}_{\tau \in \mathcal{T}(t,\ldots,T)}
        E\big[f_{\tau}(X_{\tau}) \mid \mathcal{F}_t\big],
\end{equation}
where the supremum is taken over the family
$\mathcal{T}(t,\ldots,T)$ of all $\mathbb{F}$-stopping times with
values in ${t,\ldots,T}$. Adding a positive constant $\eps$ to the
payoff $f$ just increases $V_t$ by $\eps$. We therefore can assume
without loss of generality that $f \in L_1(\mathbf{X})$ is a
positive payoff function. A stopping rule $\tau_t^* \in
\mathcal{T}(t,\ldots,T)$ is optimal for time $t$ if it attains the
optimal value
\begin{equation}\label{equ:opt stop rule}
    V_t = E\big[f_{\tau_t^*}(X_{\tau_t^*}) \mid \mathcal{F}_t\big].
\end{equation}
Once the value process is known, an optimal stopping rule at time
$t$ is given by
\begin{equation}\label{equ:opt stop times}
    \tau_t^* = \inf\{s \geq t \mid V_s \leq f_s(X_s)\}.
\end{equation}
To exploit the Markov property of the underlying process $X_t$ we
introduce the value function
\begin{equation}\label{equ:opt stop problem 2}
    v_t(x) = {\sup}_{\tau \in \mathcal{T}(t,\ldots,T)}
        E \big[f_{\tau}(X_{\tau}) \mid X_t = x \big].
\end{equation}
The Markov property implies $V_t = v_t(X_t)$. Closely related to
the value process $V_t$ is the process
\begin{equation}\label{equ:Q function}
    Q_t = {ess \sup}_{\tau \in \mathcal{T}(t+1,\ldots,T)} E\big[ f_{\tau}(X_{\tau}) \mid \mathcal{F}_t\big]
        = E\big[ f_{\tau_{t+1}^*}(X_{\tau_{t+1}^*}) \mid \mathcal{F}_t
        \big],
\end{equation}
which is defined for all $t = 0, \ldots T-1$. Again, by the Markov
property, we get the representation $Q_t = q_t(X_t)$ where
\begin{equation}\label{equ:q function 2}
    q_t(x) = {\sup}_{\tau \in \mathcal{T}(t+1,\ldots,T)}
        E\big[f_{\tau}(X_{\tau}) \mid X_t = x \big] =
        E\big[ f_{\tau_{t+1}^*}(X_{\tau_{t+1}^*}) \mid X_t = x\big].
\end{equation}
We extend the definition of $q_t$ up to the horizon $T$ and set
$q_T = f_T$. The function $q_t$ is referred to as the continuation
value. It represents the optimal value at time $t$, subject to the
constraint of not stopping at $t$. The value function and the
continuation value are related by
\begin{equation}\label{equ:Q from V}
    v_t(X_t) = \max\big(f_t(X_t), q_t(X_t)\big), \quad
    q_t(X_t) = E\big[v_{t+1}(X_{t+1}) \mid X_t \big].
\end{equation}
The dynamic programming principle implies a recursive expression
for the value, the continuation value, and the optimal stopping
times. The recursion starts at the horizon $T$ with $v_T(X_T) =
q_T(X_T) = f_T(X_T)$ and proceeds backward for $t=T-1, \ldots, 0$
according to
\begin{equation}\label{equ:V backward}
    v_t(X_t) = \max\big(f_t(X_t), E[ v_{t+1}(X_{t+1}) \mid
    X_t]\big),
\end{equation}
respectively
\begin{equation}\label{equ:Q backward}
    q_t(X_t) = E\big[\max(f_{t+1}(X_{t+1}), q_{t+1}(X_{t+1})) \mid
    X_t\big].
\end{equation}
Similarly, the recursion for the optimal stopping rules $\tau_t^*$
starts at the horizon $T$ with $\tau_T^* = T$. Given $v_t$
respectively $q_t$ and the optimal stopping rule $\tau_{t+1}^*$ at
time $t+1$, the optimal stopping rule $\tau_t^*$ is determined by
\begin{eqnarray}\label{equ:stop rule backward}
    \tau_t^* &=& t \, 1_{\{v_t(X_t) = f_t(X_t)\}} + \tau_{t+1}^*  1_{\{v_t(X_t) > f_t(X_t)\}} \nonumber \\[2mm]
             &=& t \, 1_{\{q_t(X_t) \leq f_t(X_t)\}} + \tau_{t+1}^*  1_{\{q_t(X_t) > f_t(X_t)\}}.
\end{eqnarray}
From a theoretical point of view, the value function $v_t$ and the
continuation value $q_t$ are equivalent since they both provide a
solution to the optimal stopping problem. However, from an
algorithmic point of view, the continuation value is preferred.
Indeed, $q_t$ tends to be smoother than $v_t$ because the $\max$
operation introduces a kink in the value function. We note that in
continuous time this kink disappears, since by the smooth fit
principle, the value function connects $C^1$-smoothly to the
payoff function along the optimal stopping boundary.

Expression \eqref{equ:stop rule backward} for the optimal stopping
rule suggests that we consider stopping rules parameterized by
functions $h \in L_1(\mathbf{X})$ with $h_T = f_T$. The terminal
condition $h_T = f_T$ reflects the terminal boundary condition
$\tau_T^* = T$. Let
\begin{equation}\label{equ:Markov stop time theta}
    \theta_{f,t}(h) = \theta(f_t - h_t), \quad \theta^-_{f,t}(h) = 1 - \theta(f_t - h_t),
\end{equation}
where $\theta(s) = 1_{\{s \geq 0\}}$ is the heaviside function.
Set $\tau_T(h) = T$ and define recursively
\begin{equation}\label{equ:Markov stop time}
    \tau_t(h)(\mathbf{x}) = t \, \theta_{f,t}(h)(x_t) +
        \tau_{t+1}(h)(\mathbf{x}) \, \theta^-_{f,t}(h)(x_t),
        \quad \mathbf{x} \in \bs{\mathcal{X}}.
\end{equation}
For every $h \in L_1(\mathbf{X})$ we get a valid stopping rule
$\tau_t(h)$ which does not anticipate the future, because at each
point in time $t$, the knowledge of $X_t$ is sufficient to decide
whether to stop or to continue.

\begin{defn}
The family of stopping rule $\{\tau_t(h) \mid h \in
L_1(\mathbf{X}), \; h_T = f_T\}$ is called the set of \emph{Markov
stopping rules}.
\end{defn}

The stopping rule $\tau_{t}(h)$ depends only on $h_t, \ldots,
h_{T-1}$ and is therefore constant as a function of the arguments
$x_0, \ldots, x_{t-1}$. Moreover, the recursion formula
\eqref{equ:stop rule backward} implies that the optimal stopping
rule $\tau_t^*$ at time $t$ is identical to the Markov stopping
rule $\tau_t(q)$.

\noindent Applying the Markov stopping rule $\tau_t(h)$ leads to
the cash flow $f_{\tau_{t}(h)} (X_{\tau_{t}(h)})$. More generally,
we define for $\mathbf{x} \in \bs{\mathcal{X}}$, any $0 \leq w
\leq T-t$, and $h \in L_1(\mathbf{X})$ with $h_T = f_T$ the
function
\begin{equation}\label{equ:z ipol}
    \vartheta_{t:w}(f, h)(\mathbf{x}) =
        \sum_{s=t}^{t+w} f_s(x_s) \, \theta_{f,s}(h)(x_s) \prod_{r=t}^{s-1} \theta^-_{f,r}(h)(x_r) +
             h_{t+w}(x_{t+w}) \prod_{r=t}^{t+w}
             \theta^-_{f,r}(h)(x_r),
\end{equation}
where we follow the convention that the product over an empty
index set is equal to one. The function $\vartheta_{t:w}(f, h)$
has a natural financial interpretation. It is the cash flow we
would obtain by holding the American option for at most $w$
periods, applying the stopping rule $\tau_t(h)$, and selling the
option at time $t+w$ for the price of $h_{t+w}(X_{t+w})$, if it is
not exercised before. We call $\vartheta_{t:w}(f, h)$ the cash
flow function induced by $h$.

Equations \eqref{equ:q function 2} and \eqref{equ:Q backward}
provide two different representations of $q_t$. In terms of
$\vartheta_{t:w}(f, h)$ they can be reexpressed as follows.
Because $f_{\tau_{t+1}^*} (X_{\tau_{t+1}^*}) = f_{\tau_{t+1}(q)}
(X_{\tau_{t+1}(q)}) = \vartheta_{t+1:T-t-1}(f,q)$, \eqref{equ:q
function 2} becomes
\begin{equation}\label{equ:q estimator 1}
    q_t(X_t) = E\big[\vartheta_{t+1:T-t-1}(f,q) \mid X_t \big],
\end{equation}
whereas $\vartheta_{t+1:0}(f,q) = \max(f_{t+1}, q_{t+1})$ turns
\eqref{equ:Q backward} into
\begin{equation}\label{equ:q estimator 2}
    q_t(X_t) = E\big[\vartheta_{t+1:0}(f,q) \mid X_t \big].
\end{equation}
In fact, there is a whole family of representations, parameterized
by $w \in \{0, \ldots, T - t - 1\}$. Recursively expanding
$q_{t+1}, \ldots, q_{t+w}$ in \eqref{equ:Q backward} and using the
Markov property we find that
\begin{equation}\label{equ:q estimator 3}
    q_t(X_t) = E\big[\vartheta_{t+1:w}(f,q) \mid X_t \big],
\end{equation}
for any $0 \leq w \leq T - t - 1$.

\section{Optimal Stopping as a Recursive Statistical Learning Problem}
\label{sec:opt stop stat learning}

The calculation of the recursive series of nested regression
problems \eqref{equ:q estimator 3} is becoming increasingly
demanding for high dimensional state spaces. A further
complication is introduced if the transition densities of the
Markov process $\mathbf{X}$ are not explicitly available. In this
case, the only means to assess the distribution of the Markov
process is by simulating a large number of independent sample
paths $\mathbf{X}_1, \mathbf{X}_2, \ldots, \mathbf{X}_n$. These
kind of problems are considered in statistical learning theory.

\subsection{Dynamic Look-Ahead Algorithm}

Assume a payoff $f \in  L_2(\mathbf{X})$. We interpret the unknown
continuation value $q_t \in L_2(\mathbb{R}^m, \mu_t)$ as an
approximation of the unknown optimal cash flow
$\vartheta_{t+1:w}(f,q)$, in the sense that it only depends on the
state of the underlying Markov process at time $t$. To reduce the
problem further we choose for every $t \geq 0$ a suitable set of
functions $\mathcal{H}_t$ defined on $\mathbb{R}^m$. Let
\begin{equation}\label{equ:hyp space}
    \mathcal{H} = \{h = (h_0,\ldots,h_T) : \bs{\mathcal{X}} \rightarrow \mathbb{R}^{T+1} \mid h_t \in \mathcal{H}_t \}.
\end{equation}
Given a finite amount of independent sample paths
\begin{equation}\label{equ:data}
    D_n = \{\mathbf{X}_1, \ldots, \mathbf{X}_n\},
\end{equation}
we want to find a learning rule $\hat{q}_{\mathcal H}$, i.e., a
map
\begin{equation}\label{equ:learning rule}
    \hat{q}_{\mathcal H}: D_n \mapsto \hat{q}_{\mathcal H}(D_n) =
    (\hat{q}_{\mathcal H, 0}(D_n), \ldots, \hat{q}_{\mathcal H,
    T}(D_n)) \in \mathcal{H},
\end{equation}
such that $\hat{q}_{\mathcal H, t}(D_n)$ provides an accurate
approximation of $\vartheta_{t+1:w}(f,q)$ in $\mathcal{H}_t$. The
dynamic programming principle imposes consistency conditions on a
learning rule.
\begin{defn}
A learning rule $\hat{q}_{\mathcal H}$ is called admissible if
$\hat{q}_{\mathcal H, T}(D_n) \equiv f_T$ and $\hat{q}_{\mathcal
H, t}(D_n)$, as a function of $D_n$, does not depend on the sample
paths up to and including time $t-1$, or equivalently, is a
function of $\{X_{i,s} \mid s \geq t, i=1,\ldots,n\}$ alone.
\end{defn}
We apply empirical risk minimization to recursively define an
admissible learning rule as follows. At the horizon $T$ we set
\begin{equation}\label{equ:consistency at T}
    \hat{q}_{\mathcal H, T}(D_n) \equiv f_T.
\end{equation}
For $t < T$, equation \eqref{equ:q estimator 3} suggests that we
approximate the cash flow function
\begin{equation}\label{equ:dynamic step look-ahead}
    \vartheta_{t+1:w}(f, \hat{q}_{\mathcal{H}}(D_n)),
\end{equation}
for some suitably selected parameter $w = w(t) \in \{0, \ldots,
T-t-1\}$. We choose
\begin{equation}\label{equ:ERM}
\begin{split}
     \hat{q}_{\mathcal H, t}(D_n) &= \argmin_{h \in \mathcal{H}_{t}} P_n |h - \vartheta_{t+1:w}(f, \hat{q}_{\mathcal{H}}(D_n))|^2 \\
        & = \argmin_{h \in \mathcal{H}_{t}}
            \frac{1}{n} \sum_{i=1}^n |h(X_{i,t}) - \vartheta_{t+1:w}(f,
            \hat{q}_{\mathcal{H}}(D_n))(\mathbf{X}_i)|^2,
\end{split}
\end{equation}
which is an element of $\mathcal{H}_{t}$ with minimal empirical
$L_2$-distance from the cash flow function \eqref{equ:dynamic step
look-ahead}. Because the objective function in the optimization
problem \eqref{equ:ERM} depends solely on the functions
$\hat{q}_{\mathcal H, s}(D_n)$, $s = t+1, \ldots, t+w+1$, we see
by induction that the empirical risk minimization algorithm
\eqref{equ:ERM} indeed leads to an admissible learning rule.

\begin{rem}\label{rem:dependency on Dn}
It is important to note that, while the function
$\hat{q}_{\mathcal H}(D_n)$ is a function of $\mathbf{x} \in
\bs{\mathcal{X}}$, its choice depends on the sample $D_n$.
Therefore, $\hat{q}_{\mathcal H}(D_n)$ is a random element with
values in $\mathcal{H}$ which is defined on the countable product
space $(\bs{\mathcal{X}}^{\infty}, P^{\infty},
\mathcal{F}^{\infty})$. Strictly speaking, for a sample size $n$
only the first $n$ coordinates of $\bs{\mathcal{X}}^{\infty}$ are
relevant. Analogously, the expectation
\begin{equation}\label{equ:rem dependency on Dn}
    E[\hat{q}_{\mathcal H}(D_n)] = \int_{\bs{\mathcal{X}}} \hat{q}_{\mathcal
    H}(D_n)(\mathbf{x})dP(\mathbf{x})
\end{equation}
of $\hat{q}_{\mathcal H}(D_n)$ over the path space
$\bs{\mathcal{X}}$ is still a random variable on
$\bs{\mathcal{X}}^{\infty}$.
\end{rem}

\begin{defn}\label{def:dynamic lookahead}
The dynamic look-ahead algorithm with look-ahead parameter $w =
w(t)$, $ 0 \leq w(t) \leq T-t-1$, approximates the continuation
value $q_t$ by the empirical minimizer $\hat{q}_{\mathcal H,
t}(D_n)$ of \eqref{equ:ERM}.
\end{defn}
The cash flow \eqref{equ:dynamic step look-ahead} depends on the
next $w+1$ time periods, hence, it ``looks ahead'' $w+1$ periods.
The algorithm is called ``dynamic'' because the look-ahead
parameter $w$ may be chosen time and sample dependent. We simplify
our notation and drop the explicit dependency on the sample $D_n$,
the sample size $n$, and the look-ahead parameter $w$, writing
$\hat{q}_{\mathcal H, t}$ for the solution of the empirical
minimization problem \eqref{equ:ERM}.

\subsection{Tsitsiklis-Van Roy and Longstaff-Schwartz Algorithm}
\label{sec:comparis TS and LS}

Both the Tsitsiklis-Van Roy and the Longstaff-Schwartz algorithm
are special instances of the dynamic look-ahead algorithm. The
Longstaff-Schwartz algorithm is based on the cash flow function
\begin{equation}\label{equ:target LS}
    \vartheta_{t+1}^{LS} =
    f_{\tau_{t+1}(\hat{q}_{\mathcal{H}})}\left(X_{\tau_{t+1}(\hat{q}_{\mathcal{H}})}\right),
\end{equation}
which corresponds to the maximal possible value $w = T-t-1$. On
the other extreme, the choice $w=0$ in \eqref{equ:dynamic step
look-ahead} results in the much simpler expression
\begin{equation}\label{equ:target van Roy}
    \vartheta_{t+1}^{TR} = \max(f_{t+1}, \hat{q}_{\mathcal{H},t+1}),
\end{equation}
used in the Tsitsiklis-Van Roy algorithm. In its initial form,
this algorithm has been developed to solve infinite horizon
optimal stopping problems of ergodic Markov processes. The
advantage of $\vartheta_{t+1}^{TS}$ is its numerical simplicity.
On the other hand, $\vartheta_{t+1}^{LS}$ is better suited to
approximate the optimal stopping rule because it incorporates all
future time points up to the final horizon. This property is
particularly important for Markov process with slow mixing
properties.

The dynamic look-ahead algorithm introduced in Definition
\ref{def:dynamic lookahead} interpolates between the
Tsitsiklis-Van Roy and the Longstaff-Schwartz algorithm. A dynamic
adjustment of the look-ahead parameter $w = w(t)$ allows us to
combine the algorithmic simplicity of Tsitsiklis-Van Roy and the
good approximation properties of the Longstaff-Schwartz approach.
For instance we may increase $w(t)$ for the last few time steps to
compensate the slow mixing of the Markov process.

\section{Main Results}
\label{sec:main results}

In our definition of the dynamic look-ahead algorithm
\eqref{equ:ERM} we did not further specify the approximation
scheme. The richer the set of functions $\mathcal{H}_t$, the
better it can approximate the optimal cash flow. On the other hand
large sets $\mathcal{H}_t$ would require an abundance of samples
to get a minimizer in \eqref{equ:ERM} with reasonably small
variance. These conflicting objectives are generally referred to
as the bias-variance trade-off. To get a reasonable convergence
behavior of the dynamic look-ahead algorithm, we need to impose
some restrictions on the massiveness of the approximation spaces
$\mathcal{H}_t$ and relate it to the number of samples which are
used to calculate the minimizers in \eqref{equ:ERM}.

The massiveness of a set of functions can be measured in terms of
covering and entropy numbers. The calculation of covering numbers
of classes of function has a long history dating back to
\citeasnoun{Kolmogorov-Tikhomirov:1959} and
\citeasnoun{Birman-Solomyak:1967}. We refer to
\citeasnoun{Carl-Stephani:1990} for a modern approach and
additional references. An important type of function classes for
which covering numbers can be estimated with combinatorial
techniques are the so called Vapnik-Chervonenkis classes or
VC-classes, which are by definition classes of functions of finite
VC-dimension. Informally speaking, the VC-dimension measures the
size of nonlinear sets of functions by looking at the maximum
number of sign alternations of its elements. To give a precise
definition we consider a class of functions $\mathcal{G}$ defined
on some set $S$. A set of $n$ points $\{x_1,\ldots, x_n\} \subset
S$ is said to be shattered by $\mathcal{G}$ if there exists $r \in
\mathbb{R}^n$ such that for every $b \in \{0,1\}^n$, there is a
function $g \in \mathcal{G}$ such that for each $i$, $g(x_i)
> r_i$ if $b_i = 1$, and $g(x_i) \leq r_i$ if $b_i = 0$. The
VC-dimension $\vc(\mathcal{G})$ of $\mathcal{G}$ is defined as the
cardinality of the largest set of points which can be shattered by
$\mathcal{G}$. The function classes that will appear in the
analysis of the fluctuations of the empirical minimizers
\eqref{equ:ERM} very well fit in the theory of
Vapnik-Chervonenkis. We introduce the necessary tools of the
VC-theory on the way as we prove the main results in Section
\ref{sec:proofs}.

Our error decomposition crucially depends on the convexity and the
uniform boundedness of the class of functions $\mathcal{H}_t$. We
will impose for all $t \geq 0$ the following three conditions.
\begin{enumerate}
    \item[$(\mathbf{H}_1)$] The class $\mathcal{H}_t$ is a closed convex subset of
          $L_p(\mathbb{R}^m, \mu_t)$ for some $2 \leq p \leq
          \infty$.
    \item[$(\mathbf{H}_2)$] There exists a constant $d$ such that the VC-dimension of
          $\mathcal{H}_t$ satisfies $\vc(\mathcal{H}_t) \leq d <
          \infty$.
    \item[$(\mathbf{H}_3)$] The class $\mathcal{H}_t$ is uniformly
    bounded, i.e., for some constant $H$, $|h_t| \leq H < \infty ~
          \forall h_t \in \mathcal{H}_t$.
\end{enumerate}
The convexity and uniform boundedness assumptions $(\mathbf{H}_1)$
respectively $(\mathbf{H}_3)$ are somewhat restrictive, but
encompass many common approximation schemes such as bounded convex
sets in finite dimensional linear spaces, local polynomial
approximations, or tensor product splines.

\subsection{Consistency and Convergence}
\label{sec:consist and conv}

The payoff function of an optimal stopping problem is often
unbounded. For example, in option pricing even the simplest payoff
functions of American put and call options increase linearly in
the underlying. On the other hand, any numerical algorithm works
at finite precision and tight error or convergence rate estimates
rely on some sort of boundedness assumptions. We therefore
introduce the truncation operator $T_{\beta}$ which assigns to a
real valued functions $g$ the bounded function
\begin{equation}\label{equ:trunc op}
    T_{\beta} g =
    \left \{
    \begin{array}{ll}
      g,              & \text{if } |g| \leq \beta,  \\[1mm]
      \sign(g) \beta, & \text{else}, \\
    \end{array}
    \right.
\end{equation}
and to $g \in L_p(\mathbf{X})$ its coordinate-wise truncation
$T_{\beta}g = (T_{\beta}g_0,\ldots,T_{\beta}g_T)$. We then replace
the estimator \eqref{equ:ERM} by
\begin{equation}\label{equ:ERM trunc}
    \hat{q}_{\mathcal{H}_n,t} = \hat{q}_{\mathcal{H}_n,t}(D_n) =
    \argmin_{h \in \mathcal{H}_{n,t}} P_n |h - \vartheta_{t+1:w(t)}(T_{\beta_n}f,
    \hat{q}_{\mathcal{H}_n}(D_n))|^2,
\end{equation}
where $T_{\beta_n}f$ is the payoff truncated at a threshold
$\beta_n$. The estimator \eqref{equ:ERM trunc} rests on the
hypothesis that whenever $\hat{q}_{\mathcal H_n, s}(D_n)$ is an
approximation of $q_s$ for $s \geq t+1$, then the cash flow
$\vartheta_{t+1:w}(T_{\beta_n}f, \hat{q}_{\mathcal{H}_n}(D_n))$ is
a sufficiently accurate substitute for the unknown optimal cash
flow $\vartheta_{t+1:w}(T_{\beta_n}f,q)$. We justify this
hypothesis in Proposition \ref{prop:continuity z functional} by
proving a conditional Lipschitz continuity of the functional $h
\mapsto \vartheta_{t+1:w}(T_{\beta_n}f, h)$ at $q$. The error
propagation of the recursive estimation procedure is resolved in
Corollary \ref{cor:first err decomp}, which relies on the
convexity of the approximation architecture.

The first main result provides a sufficient condition on the
growth of the number of sample paths $n$, the VC-dimension
$\vc(\mathcal{H}_{n,t})$ of the approximation spaces
$\mathcal{H}_{n,t}$, and the truncation level $\beta_n$ to ensure
convergence. Let $(\bs{\mathcal{X}}^{\infty}, P^{\infty},
\mathcal{F}^{\infty})$ be the countable product space introduced
in Remark \ref{rem:dependency on Dn}. We use the notation
$\mathbb{P} = P^{\infty}$ and denote by $\mathbb{E}$ the
expectation with respect to $\mathbb{P}$.

\begin{thm}\label{thm:ERM consistency}
Assume the payoff $f$ is in $L_2(\mathbf{X})$ and $\mathcal{H}_n$
is a sequence of approximation spaces uniformly bounded by
$\beta_n$ such that $\cup_{n=1}^{\infty} \mathcal{H}_n$ is dense
in $L_2(\mathbf{X})$. Furthermore, assume that each
$\mathcal{H}_{n,t}$ is closed, convex, and $\vc(\mathcal{H}_{n,t})
\leq d_n$. Let $\hat{q}_{\mathcal{H}_n,t}$ be the empirical
$L_2$-minimizer from \eqref{equ:ERM trunc} for a look-ahead
parameter $0 \leq w(t) \leq T-t-1$. Under the assumptions
\begin{equation}\label{equ:ERM consistency conv in prob cond}
    \beta_n \rightarrow \infty, \quad d_n \rightarrow \infty, \quad
    \frac{d_n \beta_n^2 \log(\beta_n)}{n} \rightarrow 0 \quad (n \rightarrow
    \infty),
\end{equation}
it follows that
\begin{equation}\label{equ:ERM consistency conv in prob}
    \|\hat{q}_{\mathcal{H}_n,t} - q_t \|_2 \rightarrow 0,
\end{equation}
in probability and in $L_1(\mathbb{P})$. If furthermore
\begin{equation}\label{equ:ERM consistency conv almost sure cond}
    \frac{\beta_n^2 \log(n)}{n} \rightarrow 0,
\end{equation}
then the convergence in \eqref{equ:ERM consistency conv in prob}
holds almost surely.
\end{thm}

\begin{proof}
See Section \ref{sec:proof consistency thm}.
\end{proof}

Theorem \ref{thm:ERM consistency} proves convergence of the
truncated version \eqref{equ:ERM trunc} of the dynamic look-ahead
algorithm. It generalizes previous results in two directions.
First, the number of samples, the size of the approximation
architecture (measured in terms of the VC-dimension), and the
truncation threshold are increased simultaneously.
\citeasnoun{Glass-Yu:2003} address the same question for the
Longstaff-Schwartz algorithm with linear finite dimensional
approximation. They avoid truncation by imposing fourth-order
moment conditions and find that the number samples must grow
surprisingly fast. For example, if $X_t$ is log-normally
distributed and $n$ denotes the dimension of the linear
approximation space, the number of samples must be proportional to
$\exp(n^2)$. Second, Theorem \ref{thm:ERM consistency} covers
approximation architectures of bounded VC-dimension and does not
depend on the law of the underlying Markov process. For instance,
the convergence proof of \citeasnoun{Protter-Lamberton:2001}
relies on the additional assumption $P(q = f) = 0$.

In \eqref{equ:ERM trunc} we reduce unbounded to bounded payoffs by
truncating at a suitable cutoff level. The next result bounds the
approximation error in terms of the cutoff level.
\begin{prop}\label{prop:truncation}
Let $1 \leq p < \infty$ and $f \in L_p(\mathbf{X})$ be a
nonnegative payoff function. If $\bar{q}_{\beta}$ is the
continuation value of the truncated payoff~ $T_{\beta} f$, it
follows that
\begin{equation}\label{equ:truncation conv}
    \| q_t - \bar{q}_{\beta,t} \|_p \rightarrow 0,
\end{equation}
for $\beta \rightarrow \infty$, and if $1 < r < p$, then
\begin{equation}\label{equ:truncation rate}
    \| q_t - \bar{q}_{\beta,t} \|_r \leq
        \sum_{s=t+1}^{T} \Big( r \int_{\beta}^{\infty} u^{r-1} P(f_{t+1} > u) \, d u \Big)^{\frac{1}{r}} \leq O(\beta^{\frac{r-p}{r}}).
\end{equation}
\end{prop}

\begin{proof}
See Section \ref{sec: proof trunc error estimate}.
\end{proof}

The bound \eqref{equ:truncation rate} can be refined in terms of
Orlicz norms. The Orlicz norm of a random variable $Y$ is defined
as
\begin{equation}\label{equ:Orlicz norm}
    \| Y \|_{\psi} = \inf \{ C > 0 \mid  E \left[\psi\left(|Y|C^{-1}\right)\right] \leq  1 \},
\end{equation}
where $\psi$ is a nondecreasing, convex function with $\psi(0) =
0$. Note that $\psi(y) = y^p$ reduces to the usual $L_p$-norms. If
$\| f_{t+1} \|_{\psi} < \infty$ Markov's inequality implies the
tail bound
\begin{equation}\label{equ:tail bound Orlicz}
    P(f_{t+1} > u) \leq \frac{1}{\psi(u \| f_{t+1}
    \|_{\psi}^{-1})},
\end{equation}
which we then can apply to the middle term in
\eqref{equ:truncation rate}. In particular, $\psi_p(x) = \exp(x^p)
- 1$ leads to the exponential bound
\begin{equation}\label{equ:tail bound Orlicz spec}
    P(f_{t+1} > u) \leq \exp\left(- u^p \| f_{t+1} \|_{\psi_p}^{-1}\right)
        \left(1-\exp\left(- \beta^p \| f_{t+1} \|_{\psi_p}^{-1}\right)\right)^{-1},
\end{equation}
for all $u \geq \beta$. In financial applications a typical
situation is $f_{t+1} = f(\exp(X_{t+1}))$, where $X_{t+1}$ is
normally distributed and $f(y) \leq C y^q$ has polynomial growth.
The tail estimate
\begin{equation}\label{equ:tail bound norm result}
    P(f_{t+1} > u) \leq O\left( \frac{1}{\log(u)} \exp\left(-
    \log(u)^2\right)\right)
\end{equation}
is a direct consequence of the well-known asymptotic expansion
\begin{equation}\label{equ:tail bound norm}
    1-\Phi(u) \leq \phi(u) u^{-1}\left(1 - \frac{1}{u^2} + \frac{3}{u^4} + O(u^{-6})\right)
\end{equation}
for the tail of the standard normal distribution $\Phi$ with
density $\phi$. \eqref{equ:tail bound norm result} improves the
rate of order $O(\beta^{1 - p / r})$ in \eqref{equ:truncation
rate} considerably, despite of the logarithmic terms in the
exponent.

\subsection{Error Estimate and Sample Complexity}
\label{sec:reates and complexity}

Theorem \ref{thm:ERM consistency} shows that simultaneously
increasing the truncation threshold, the VC-dimension of the
approximation architecture, and the number of samples at a proper
rate, the resulting estimator \eqref{equ:ERM trunc} converges to
the solution of the optimal stopping problem. Proposition
\ref{prop:truncation} quantifies the error of an initial
truncation at a fixed threshold. We continue the error analysis of
the dynamic look-ahead algorithm by truncating unbounded payoffs
at a sufficiently large threshold $\Theta$ and considering a
single approximation architecture $\mathcal{H}$. The second main
result bounds the overall error for bounded payoff functions in
terms of the approximation error and the sample error,
generalizing the familiar bias-variance trade-off in nonparametric
regression and density estimation.

\begin{thm}\label{thm:ERM conv rate}
Consider a payoff $f \in L_{\infty}(\mathbf{X})$ with
$\|f_t\|_{\infty} \leq \Theta$. Assume that each $\mathcal{H}_{t}$
is a closed convex set of functions, uniformly bounded by $H$,
with $\vc(\mathcal{H}_{t}) \leq d$. Let
$\hat{q}_{\mathcal{H},t}(D_n)$ be the empirical $L_2$-minimizer
from \eqref{equ:ERM} for a look-ahead parameter $0 \leq w(t) \leq
T-t-1$. Set $\beta = max(\Theta, H)$. Then, for $n \geq
382\beta^2/\eps$,
\begin{eqnarray}\label{equ:ERM conv rate error estimate}
   \mathbb{E} \left[\|\hat{q}_{\mathcal{H},t}(D_n) - q_t \|_2^2 \right] &\leq&
          2 \cdot 16^{w(t)} \max_{s=t,\ldots,t+w(t)+1} \; \inf_{h \in \mathcal{H}_{s}} \|h - q_s \|_2^2 \; +  \\[2mm]
          &&  2 \cdot 16^{w(t)} (w(t)+2) \left(\frac{6998 \beta^2 + \log(6998 K \beta^2)}{n} + \frac{v
                \log(n)}{n}\right), \nonumber
\end{eqnarray}
where
\begin{equation*}\label{equ:ERM conv rate error estimate 1}
    v = 2d(c(w(t)) + 1), \quad K = 6 e^4(d+1)^2(c(w(t))d+1)^2(1024 e \beta)^{v},
\end{equation*}
and
\begin{equation*}\label{equ:ERM conv rate error estimate 2}
    c(w(t)) = 2(w(t)+2)\log_2(e(w(t)+2)).
\end{equation*}
\end{thm}
\begin{proof} See Section \ref{sec:proof consistency thm}. \end{proof}

The effectiveness of a learning algorithm can be quantified by the
number of samples which are required to produce with high
confidence $1 - \delta$ an almost minimizer
\begin{equation}\label{equ:almost min samp complexity}
    \|\hat{q}_{\mathcal{H},t}(D_n) - q_t \|_2^2 \leq
        \inf_{h_t \in \mathcal{H}_{t}} \|h_t - q_t \|_2^2 + \eps, \quad \forall \,
        t=0\ldots,T-1,
\end{equation}
for a certain error accuracy $\eps$. In \eqref{equ:almost min samp
complexity} the error is measured relative to the minimal
approximation error at time step $t$. It is evident from
\eqref{equ:ERM conv rate error estimate} that an accurate estimate
is only obtained if the approximation error in all previous
learning tasks is small as well. To disentangle sample complexity
and approximation error, we measure the performance of the
learning rule relative to the overall approximation error in
\eqref{equ:ERM conv rate error estimate}.

\begin{cor}\label{cor:sample complexity bound}
Assume $f \in L_{\infty}(\mathbf{X})$ with $\|f_t\|_{\infty} \leq
\Theta$ and let $\mathcal{H}$ be as in Theorem \ref{thm:ERM conv
rate}. The sample complexity
\begin{equation}\label{equ:sample complexity bound}
\begin{split}
    c(\eps, \delta) = \min & \Bigg\{ \; n_0 \; \big | \; \forall n \geq n_0, \\
        &\mathbb{P} \left( \|\hat{q}_{\mathcal{H},t}(D_n) - q_t \|_2^2
            \geq 2 \cdot 16^{w(t)} \max_{s=t,\ldots,t+w(t)+1} \; \inf_{h \in \mathcal{H}_{s}} \|h - q_s \|_2^2 +
            \eps \right) \leq \delta \Bigg\}
\end{split}
\end{equation}
of the empirical $L_2$-minimizer \eqref{equ:ERM} is bounded by
\begin{equation}\label{equ:sample complexity bound 1}
    c(\eps, \delta) \leq 2 \cdot 13996 (w(t)+2)16^{w(t)}
    \beta^2 \max\left(\frac{1}{\eps} \log\left( \frac{K}{\delta}\right), v \log\left( \frac{1}{\eps}\right)
    \right),
\end{equation}
where $\beta$, $v$ and $K$ are as in Theorem \ref{thm:ERM conv
rate}.
\end{cor}

\begin{proof} See Section \ref{sec:proof consistency thm}. \end{proof}

Theorem \ref{thm:ERM conv rate} and Corollary \ref{cor:sample
complexity bound} estimate the sample error for a fixed
approximation scheme and truncation threshold. The bound
\eqref{equ:ERM conv rate error estimate} and the complexity
estimate \eqref{equ:sample complexity bound 1} hold uniformly for
any law of $\mathbf{X}$ and payoff function $f$ with
$\|f\|_{\infty} \leq \Theta$. Hence, the bounds are independent of
the distribution of the underlying Markov process, the optimal
stopping time, and the smoothness of the continuation value. The
asymptotic rate $O(\log(n) n^{-1})$ of the sample error (the
second term on the right hand side of \eqref{equ:ERM conv rate
error estimate}) is typical for nonparametric least square
estimates with approximation schemes of finite VC-dimension, see,
e.g., \citeasnoun[Theorem 11.5]{Kohler-etal:2002}.

If we impose additional assumptions on the smoothness of the
continuation value $q$ the approximation errors $\inf_{h \in
\mathcal{H}_{n,s}} \|h - q_s \|_2^2$ in \eqref{equ:ERM conv rate
error estimate} can be estimated further by approximation theory.
Smoothness assumptions are not unreasonable. Although for many
financial applications the payoff is only continuous or piecewise
continuous, the continuation value is often smooth. The degree of
smoothness of $q$ is crucial for how to choose approximation
spaces $\mathcal{H}_n$ to get the most favorable rate of
convergence by properly balancing the approximation error and the
sample error.

Smoothness is often measured in terms of Sobolev spaces
$W^k(L_{p}(\Omega, \lambda))$, where $\Omega \subset \mathbb{R}^m$
is a domain in $\mathbb{R}^m$ and $\lambda$ is the Lebesgue
measure on $\Omega$. These are functions $g \in L_{p}(\Omega,
\lambda)$ which have all their distributional derivatives of order
up to $k$ in $L_{p}(\Omega, \lambda)$. The Sobolev (semi-)norm $\|
g \|_{p,k,\Omega,\lambda}$ may be regarded as a measure of
smoothness for a function $g \in W^k(L_{p}(\Omega, \lambda))$.

In practical applications of the Longstaff-Schwartz algorithm
approximation by polynomials performs rather well. Let
$\mathcal{P}_r$ be the space of multivariate polynomials on
$\mathbb{R}^m$ with coordinate wise degree at most $r-1$. For
simplicity we assume $X_t$ is localized to a sufficiently large
cube $I \subset \mathbb{R}^m$. This assumption can be satisfied by
applying a truncation argument similar to the one developed in
Proposition \ref{prop:truncation}.

\begin{cor}\label{cor:Sobolev approx rates}
Assume that $X_t$ is localized to a cube $I \subset \mathbb{R}^m$,
$f \in L_{\infty}(\mathbf{X})$, and that the continuation value
$q_t$ is in the Sobolev space $W^k(L_{\infty}(I, \lambda))$ for
all $t$. Define the sequence of approximation architectures
\begin{equation}\label{equ:Sobolev approx rates approx space}
    \mathcal{H}_{n,t} = \{ p \in \mathcal{P}_{n^{1/(m+2k)}} \mid \|p\|_{\infty,I,\lambda} \leq
        2\|q_t\|_{\infty,k,I,\lambda}\}.
\end{equation}
Then,
\begin{equation}\label{equ:Sobolev approx rates cor}
    \mathbb{E} \left[\|\hat{q}_{\mathcal{H}_n,t}(D_n) - q_t \|_2^2
    \right] \leq O\left(\log(n) n^{-\frac{2k}{2k + m}}\right).
\end{equation}
If $\mu_t$ has a bounded density with respect to the Lebesgue
measure and $q_t \in W^k(L_{p}(I, \lambda))$ for some $p \geq 2$
the same result holds if we replace $\mathcal{H}_{n,t}$ in
\eqref{equ:Sobolev approx rates approx space} by
\begin{equation}\label{equ:Sobolev approx rates approx space 2}
    \mathcal{H}_{n,t} = \{ p \in \mathcal{P}_{n^{1/(m+2k)}} \mid \|p\|_{p,I,\lambda} \leq
        2\|q_t\|_{p,k,I,\lambda}\}.
\end{equation}
\end{cor}

\begin{proof} The result essentially follows from
Jackson type estimates, Theorem 6.2 in Chapter 7 of
\citeasnoun{DeVore-Lorentz:1993}. See Section \ref{sec: proof
Sobolev approx rates}.
\end{proof}

Corollary \ref{cor:Sobolev approx rates} is a prototypical
application of Theorem \ref{thm:ERM conv rate} to global
approximation by polynomials. Other approximation schemes can be
treated similarly, as long as the conditions
$(\mathbf{H}_1)$-$(\mathbf{H}_3)$ are satisfied. To get the rate
stated in Corollary \ref{cor:Sobolev approx rates} the dimension
$n^{m/(m+2k)}$ of the polynomial approximation architecture
\eqref{equ:Sobolev approx rates approx space} has to grow with
increasing sample size such that the approximation error and the
sample error are balanced. The rate \eqref{equ:Sobolev approx
rates cor} is up to a logarithmic term the lower minimax rate of
convergence for estimating regression functions, see
\citeasnoun{Stone:1982}.

\subsection{Discussion and Remarks}

The Longstaff-Schwartz algorithm and its generalization, the
dynamic look-ahead algorithm, perform surprisingly well for many
practical applications such as pricing American options which are
not too far in or out of the money. This empirical observation can
be explained as follows. It follows from \eqref{equ:q estimator 3}
that an approximation of the optimal cash flow
$\vartheta_{t+1:w}(f, q)$ can be used to estimate the continuation
value at time $t$. A closer look at definition \eqref{equ:z ipol}
shows that for the maximal possible value $w = T - t - 1$ the cash
flow $\vartheta_{t+1:w}(f, h)$ is close (in the $L_2$-sense) to
the optimal $\vartheta_{t+1:w}(f, q)$ if the signs of $f-h$ and
$f-q$ disagree only on a subset of the path space with small
probability, or equivalently if the probability of the symmetric
difference
\begin{equation}\label{equ:small sign mismatch}
    P(\{f-h > 0\} \Delta \{f-q > 0\}),
\end{equation}
is small. Note that a small probability \eqref{equ:small sign
mismatch} does not necessarily entail that the functions $h$ and
$q$ are close in the $L_2$-sense. If the look-ahead parameter $w$
satisfies $w < T - t - 1$ then $\vartheta_{t+1:w}(f, h)$ is a good
approximation of the optimal cash flow if, in addition to a small
probability \eqref{equ:small sign mismatch}, also the
$L_2$-distance between $h_{t+w+1}$ and the unknown continuation
value $q_{t+w+1}$ is small. Consequently, a look-ahead parameter
$0 \leq w < T - t - 1$ requires good approximations for
$q_{w+1},\ldots,q_{T-1}$. Determining accurate and stable
estimators for $q_t$ with $t$ close to 1 may be difficult to
achieve, in particular if the samples of the Markov process do not
cover sufficiently large parts of the state space. This explains
why the Tsitsiklis-Van Roy algorithm (corresponding to $w=0$) may
perform badly for finite horizon problems.

As opposed to the empirically demonstrated efficiency of the
Longstaff-Schwartz algorithm, the results of Theorem \ref{thm:ERM
conv rate} and Corollary \ref{cor:sample complexity bound} are
somewhat pessimistic. For practical parameter values $\eps$,
$\delta$, $d$, $w$, and large enough cutoff level $\beta$, the
sample complexity bound \eqref{equ:sample complexity bound 1}
leads to a very large sample size. The reason for the pessimistic
sample size estimates is twofold. First, the estimator
$\hat{q}_{\mathcal{H}}$ is sensitive to error propagation effects
caused by the backward induction. This leads to error estimates
such as \eqref{equ:ERM conv rate error estimate} which depend
exponentially on the number of look-ahead periods $w(t)$. The
minimal choice $w = 0$ would resolve the exponential dependence
but, as explained above, may have limited capabilities to
approximate the optimal cash flow. Another reason is the
generality of our error estimates. We already observed that
$\hat{q}_{\mathcal{H}}$ leads to an accurate approximation of the
optimal cash flow if the probability of the symmetric difference
$P(\{f-\hat{q}_{\mathcal{H}} > 0\} \Delta \{f-q > 0\})$ is small.
However, it is difficult to derive error estimates which take this
effect into account without imposing additional assumptions on the
smoothness of the payoff and the distribution of the stopping time
in the neighborhood of $\{q = f\}$.

We considered in this work estimators based on straightforward
empirical $L_2$-risk minimization. A deficiency of the simple
estimator considered in Corollary \ref{cor:Sobolev approx rates}
is that the degree of smoothness and an upper bound for
$\|q_t\|_{\infty,k,I,\lambda}$ has to be known. There exists a
variety of advanced nonparametric regression estimators which have
been developed to cope with the shortcomings of the basic
empirical risk minimization procedure. The main generalizations in
this direction are sieve estimators, studied for example by
\citeasnoun{Shen-Wong:1994}, \citeasnoun{Shen:1997}, and
\citeasnoun{Brige-Massart:1998}, adaptive methods such as
complexity regularization, penalization, and model selection, see
\citeasnoun{Barron-Brige-Massart:1999},
\citeasnoun{Kohler-etal:2002}, and the reference therein.

The benefit of conditions $(\mathbf{H}_1)$-$(\mathbf{H}_3)$ is
that convexity arguments and VC-techniques lead to error estimates
without the necessity of imposing further assumptions on the
Markov process $\mathbf{X}$. On the downside, some important
commonly used approximation schemes are excluded. For instance,
condition $(\mathbf{H}_2)$ conflicts with approximation in Sobolev
or Besov balls, which have infinite VC-dimension, and the
convexity condition $(\mathbf{H}_3)$ is incompatible with many
interesting nonlinear approximation schemes, such as $n$-term
approximation, wavelet thresholding, or neural network
architectures.

A promising approach to extend and refine the results of this work
is to approximate the cash flow $\vartheta_{t:w}(f, h)$ by a
suitably smoothed version with better Lipschitz continuity
properties. We then can express the massiveness of the
approximation schemes directly in terms of covering numbers and
exploit the dependency of the covering numbers on the radius of
the function class. The additional step of first bounding the
VC-dimension becomes unnecessary. However, this approach is of
less generality because it depends on the additional assumptions
that the probability $P(\{|q - f| < \eps \})$ decays to zero as
$\eps \rightarrow 0$ and the semi-group generated by the Markov
process $\mathbf{X}$ has good smoothing properties.

Once we have selected a sequence of approximation architectures
$\mathcal{H}_{n,t}$ the final step towards an implementation is to
determine a computationally efficient algorithm that minimizes the
empirical $L_2$-risk \eqref{equ:ERM trunc} over
$\mathcal{H}_{n,t}$ in a polynomial number of time steps.
Unfortunately, for many approximation spaces, such as certain
neural network architectures, constructing a solution which nearly
minimizes the empirical $L_2$-risk turns out to be NP-complete or
even NP-hard. Thus there might still exist serious complexity
theoretic barriers to efficient numerical implementations of
specific approximation schemes.

\subsection{Acknowledgements}

The author would like to thank Paul Glasserman, Tom Hurd, Markus
Leippold, Maung Min-Oo, and Paolo Vanini for helpful discussions.
The detailed comments and suggestions of a referee greatly helped
to improve a first version of this paper.

\section{Proofs}
\label{sec:proofs}

The proof of the main results, Theorem \ref{thm:ERM consistency}
and \ref{thm:ERM conv rate}, is divided into tree steps. The
strategy is as follows. First, we prove in Corollary
\ref{cor:first err decomp} an error decomposition in terms of an
approximation error and an expected centered loss
\eqref{equ:centered loss}. The second step is to estimate the
covering numbers of the so called centered loss class
\eqref{equ:centered loss class}, see Corollary \ref{cor:covering
number bounds}. The last step is to apply empirical process
techniques to bound the fluctuation of the expected centered loss
in terms of the covering numbers.

\subsection{Error Decomposition}

We assume from now on without further mentioning that $\mathcal{H}
\subset L_2(\mathbf{X})$ and that all approximation spaces
$\mathcal{H}_t$ are closed and convex. Before we can state our
main error decomposition we need to introduce some more notation.
Let
\begin{equation}\label{equ:proj H}
    \pi_{\mathcal{H}_t}: L_2(\mathbb{R}^d, \mu_t) \rightarrow \mathcal{H}_t
\end{equation}
denote the projection onto the closed convex subset $\mathcal{H}_t
\subset L_2(\mathbb{R}^m,\mu_t)$ and set
\begin{equation}\label{equ:extended proj}
    \pr_{\mathcal{H}_t} = \pi_{\mathcal{H}_t} \circ E[\ \cdot \mid X_t = \cdot \ ] : L_2(\bs{\mathcal{X}}, P) \rightarrow
    \mathcal{H}_t.
\end{equation}
For any $h =(h_0,\ldots,h_T) : \bs{\mathcal{X}} \rightarrow
\mathbb{R}^{T+1}$ with $h_T = f_T$ we introduce the centered loss
\begin{equation}\label{equ:centered loss}
    l_t(h) = |h_t - \vartheta_{t+1:w}(f,h)|^2 - |\pr_{\mathcal{H}_t} \vartheta_{t+1:w}(f,h) - \vartheta_{t+1:w}(f,h)|^2.
\end{equation}
In favor of a more compact notation we have dropped the dependency
of $l_t(h)$ on the look-ahead parameter $w$. Note that the
centered loss $l_t(h)$ only depends on $h_t, \ldots, h_{T-1}$ and
can take on negative values. However, $E[l_t(h)] \geq 0$ as we
will see in Lemma \ref{lem:err dec convex} below.

We decompose the overall error into an approximation error, a
sample error, and a third term which captures the error
propagation caused by the recursive definition of the dynamic
look-ahead estimator.

\begin{prop}\label{prop:first err decomp}
Assume that $\hat{q}_{\mathcal{H}}$ is the result of an admissible
learning rule. Then
\begin{eqnarray}\label{equ:first err decomp}
    \|\hat{q}_{\mathcal{H},t} - q_t \|_2 &\leq&
        \inf_{h \in \mathcal{H}_t} \|h - q_t \|_2 \ + \
        E[l_t(\hat{q}_{\mathcal{H}})]^{1/2} \ + \
        3 \ \sum_{s=t+1}^{t+w+1} \|\hat{q}_{\mathcal{H},s} - q_s \|_2.
\end{eqnarray}
\end{prop}

In general we cannot approximate
$\vartheta_{t+1:w}(f,\hat{q}_{\mathcal{H}})$ by functions $h_t \in
L_2(\mathbb{R}^d, \mu_t)$ arbitrarily well and therefore
\begin{equation}\label{equ:pos approx error}
    \inf_{h_t \in \mathcal{H}_t} E[|h_t - \vartheta_{t+1:w}(f,\hat{q}_{\mathcal{H}})|^2] > 0.
\end{equation}
For this reason we base our error decomposition \eqref{equ:first
err decomp} on the more complicated centered loss function which
expresses the sample error relative to the optimal one-step
expected loss
\begin{equation}\label{equ:one-step approx error}
    E[|\pr_{\mathcal{H}_t} \vartheta_{t+1:w}(f,\hat{q}_{\mathcal{H}})
        - \vartheta_{t+1:w}(f, \hat{q}_{\mathcal{H}}) |^2].
\end{equation}
The first term on the right-hand side of \eqref{equ:first err
decomp} is the approximation error, a deterministic quantity,
which can be analyzed by approximation theory. The second term
$E[l_t(\hat{q}_{\mathcal{H}})]^{1/2}$ is usually referred to as
the sample error. The last term in \eqref{equ:first err decomp}
collects the error propagation introduced by the previous learning
tasks through the dynamic programming backward recursion.

\begin{cor}\label{cor:first err decomp}
Let
\begin{equation}\label{equ:one step add error}
    e_t = \inf_{h \in \mathcal{H}_{t}} \|h - q_t \|_2 + E[l_t(\hat{q}_{\mathcal{H}})]^{1/2}
\end{equation}
denote the one-step error. Then,
\begin{equation}\label{equ:first err decomp recursion solved}
    \|\hat{q}_{\mathcal{H},t} - q_t \|_2 \leq e_t \ + \ 3 \sum_{s=t+1}^{t+w+1}
    4^{s-t-1}e_s,
\end{equation}
and
\begin{equation}\label{equ:first err decomp recursion solved 2}
    \|\hat{q}_{\mathcal{H},t} - q_t \|_2 \leq 4^{w+1} \max_{s=t,\ldots,t+w+1} \left(
        \inf_{h \in \mathcal{H}_{s}} \|h - q_s \|_2 + E[l_s(\hat{q}_{\mathcal{H}})]^{1/2}
        \right).
\end{equation}
\end{cor}

\begin{proof}[Proof of Corollary \ref{cor:first err decomp}]
This follows at once from \eqref{equ:first err decomp} by
recursively inserting the error estimate \eqref{equ:first err
decomp} for $s \geq t+1$.
\end{proof}

The proof of the error decomposition \eqref{equ:first err decomp}
crucially relies on the convexity of the approximation spaces,
Lemma \ref{lem:err dec convex}, and a Lipschitz estimate for
$\vartheta_{t+1:w}(f,h)$ as a function of $h$, Proposition
\ref{prop:continuity z functional}.

\begin{lem} \label{lem:err dec convex}
Denote by
\begin{equation}\label{equ:reg function}
    \rho_t(h)(x) = E[\vartheta_{t+1:w}(f,h) \mid X_t = x]
\end{equation}
the regression function of $\vartheta_{t+1:w}(f,h)$. For any $h
\in \mathcal{H}$ with $h_T = f_T$
\begin{equation}\label{equ:err decomp convex}
    \|h_t - \pi_{\mathcal{H}_t} \rho_t(h) \|_2^2 = \|h_t - \pr_{\mathcal{H}_t} \vartheta_{t+1:w}(f,h) \|_2^2
        \leq E[l_t(h)].
\end{equation}
In particular $E[l_t(h)] \geq 0$.
\end{lem}

\begin{proof} The proof is identical to the proof of Lemma 5
in \citeasnoun{Cucker-Smale:2001}. Because $\rho_t(h)$ is the
regression function of $\vartheta_{t+1:w}(f,h)$, which only
depends on $h_{t+1}, \ldots, h_{T-1}$, we have for all $h_t \in
L_2(\mathbb{R}^d, \mu_t)$
\begin{equation}\label{equ:regression func}
    \|h_t - \rho_t(h) \|_2^2 =
        E[|h_t - \vartheta_{t+1:w}(f,h) |^2 - |\rho_t(h) - \vartheta_{t+1:w}(f,h) |^2].
\end{equation}
Let $h \in \mathcal{H}$ be arbitrary. Since $\mathcal{H}_t$ is
convex and since $\pr_{\mathcal{H}_t} \vartheta_{t+1:w}(f,h) =
\pi_{\mathcal{H}_t} \rho_t(h)$ minimizes the distance to
$\rho_t(h)$ it follows that
\begin{equation}\label{equ:obtuse 1}
    \langle \rho_t(h) - \pi_{\mathcal{H}_t} \rho_t(h),
    h_t - \pi_{\mathcal{H}_t} \rho_t(h) \rangle \leq
    0.
\end{equation}
Therefore
\begin{equation}\label{equ:obtuse 2}
    \| \pr_{\mathcal{H}_t} \vartheta_{t+1:w}(f,h)  - h_t\|_2^2 = \| \pi_{\mathcal{H}_t} \rho_t(h)  - h_t\|_2^2
    \leq \| \rho_t(h) - h_t \|_2^2 -
        \| \rho_t(h) - \pi_{\mathcal{H}_t} \rho_t(h) \|_2^2.
\end{equation}
Because both $h_t$ and $\pi_{\mathcal{H}_t} \rho_t(h)$ are in
$\mathcal{H}_t$ we can apply \eqref{equ:regression func} twice
which shows that the right hand side of \eqref{equ:obtuse 2} is
equal to
\begin{equation}\label{equ:obtuse 3}
    \| \rho_t(h) - h_t \|_2^2 -
    \| \rho_t(h) - \pi_{\mathcal{H}_t} \rho_t(h) \|_2^2
    =
    E[|h_t - \vartheta_{t+1:w}(f,h) |^2 -
      |\pi_{\mathcal{H}_t} \rho_t(h) - \vartheta_{t+1:w}(f,h) |^2 ].
\end{equation}
\end{proof}

For $w = 0$ we immediately obtain from $|\max(a,x) - \max(a,y)|
\leq |x-y|$ and Jensen's inequality the uniform Lipschitz bound
\begin{equation}\label{equ:continuity TvR}
    \| E[ \vartheta_{t+1:0}(f,g) - \vartheta_{t+1:0}(f,h) \mid X_t \ ] \|_p \leq
        \| g_{t+1} - h_{t+1} \|_p.
\end{equation}
More generally, we have the following conditional Lipschitz
continuity at the continuation value.

\begin{prop}\label{prop:continuity z functional}
For every $h \in L_p(\mathbf{X})$ with $h_T = f_T$ and $0 \leq w
\leq T-t$
\begin{eqnarray}\label{equ:continuity z functional 0}
    \| E[ \vartheta_{t+1:w}(f,h) \mid X_t \ ] - q_t \|_p &=&
            \| E[ \vartheta_{t+1:w}(f,h) - \max(f_{t+1}, q_{t+1}) \mid X_t \ ] \|_p \nonumber \\
        &=& \| E[ \vartheta_{t+1:w}(f,h) - \vartheta_{t+1:w}(f,q) \mid X_t \ ] \|_p.
\end{eqnarray}
Furthermore,
\begin{equation}\label{equ:continuity z functional}
    \| E[ \vartheta_{t+1:w}(f,h) - \vartheta_{t+1:w}(f,q) \mid X_t \ ] \|_p \ \leq \
        \sum_{s=t+1}^{t+w+1}\| h_{s} - q_{s}\|_p.
\end{equation}
\end{prop}

A similar estimate for the special case $w = T-t-1$ can also be
found in \citeasnoun{Protter-Lamberton:2001}. Note that the
uniform Lipschitz estimate \eqref{equ:continuity TvR} does not
extend to $w > 0$. Proposition \ref{prop:continuity z functional}
only provides a Lipschitz estimate at the continuation value.

\begin{proof} First note that from the Markov property
\begin{equation}\label{equ:appl of Markov prop}
    E[ \vartheta_{t+1:w}(q) -\vartheta_{t+1:w}(h)\mid X_t \ ] =
    E[ \vartheta_{t+1:w}(q) - \vartheta_{t+1:w}(h)\mid \mathcal{F}_t \ ].
\end{equation}
Equation \eqref{equ:continuity z functional 0} follows directly
from the recursive definition of $q_t$. The case $w=0$ is covered
in \eqref{equ:continuity TvR}. For $w>0$ it follows from the
definition of $\vartheta_{t+1:w}$ that
\begin{multline*}
    \| E[ \vartheta_{t+1:w}(q) - \vartheta_{t+1:w}(h)\mid \mathcal{F}_t \ ] \|_p \ \leq \\
     \| E[f_{t+1} (\theta_{f,t+1}(q) - \theta_{f,t+1}(h)) +
            \theta_{f,t+1}^-(q) \vartheta_{t+2:w-1}(q) - \theta_{f,t+1}^-(h) \vartheta_{t+2:w-1}(h) \mid \mathcal{F}_t \ ] \|_p.
\end{multline*}
Adding and subtracting the term $q_{t+1}(\theta_{f,t+1}(q) -
\theta_{f,t+1}(h))$, the triangle inequality implies
\begin{multline*}
    \| E[\vartheta_{t+1:w}(q) - \vartheta_{t+1:w}(h)\mid \mathcal{F}_t \ ] \|_p
        \; \leq \; \| E[(f_{t+1} - q_{t+1}) (\theta_{f,t+1}(q) - \theta_{f,t+1}(h)) \mid \mathcal{F}_t \ ] \|_p \; + \\
        \| E[ \theta_{f,t+1}^-(q)\vartheta_{t+2:w-1}(q) - \theta_{f,t+1}^-(h)\vartheta_{t+2:w-1}(h) +
        q_{t+1}(\theta_{f,t+1}(q) - \theta_{f,t+1}(h)) \mid \mathcal{F}_t \ ] \|_p.
\end{multline*}
Now
\begin{eqnarray*}\label{equ:continuity y func functional 1}
    \theta_{f,t+1}(q) - \theta_{f,t+1}(h) &=&
         1_{\{f_{t+1} \geq q_{t+1}\}} - 1_{\{f_{t+1} \geq h_{t+1}\}} \\
     &=& 1_{\{0 \leq f_{t+1} - q_{t+1} < h_{t+1} - q_{t+1}\}} -
         1_{\{h_{t+1} - q_{t+1} \leq f_{t+1} - q_{t+1} < 0\}},
\end{eqnarray*}
which leads to
\begin{eqnarray*}
  &&\hspace{-1cm} (f_{t+1} - q_{t+1}) (1_{\{0 \leq f_{t+1} - q_{t+1} < h_{t+1} - q_{t+1}\}} -
        1_{\{h_{t+1} - q_{t+1} \leq f_{t+1} - q_{t+1} < 0\}})  \\
    &\leq& (h_{t+1} - q_{t+1}) 1_{\{h_{t+1} - q_{t+1} > 0\}} - (h_{t+1} - q_{t+1}) 1_{\{h_{t+1} - q_{t+1} < 0\}} \\
    &\leq& |h_{t+1} - q_{t+1}|.
\end{eqnarray*}
By the Markov property $q_{t+1}(X_{t+1}) = E[
\vartheta_{t+2:w-1}(q) \mid \mathcal{F}_{t+1}]$. Because
$\theta_{f,t+1}(q)$ and $\theta_{f,t+1}(h)$ are
$\sigma(X_{t+1})$-measurable it follows that
\begin{eqnarray}\label{equ:continuity y func functional 2}
     E[q_{t+1} (\theta_{f,t+1}(q) - \theta_{f,t+1}(h)) \mid \mathcal{F}_t \ ]
         &=& E[E[\vartheta_{t+2:w-1}(q) \mid \mathcal{F}_{t+1}] (\theta_{f,t+1}(q) - \theta_{f,t+1}(h)) \mid \mathcal{F}_t \ ] \nonumber \\
         &=& E[\vartheta_{t+2:w-1}(q) (\theta_{f,t+1}(q) - \theta_{f,t+1}(h)) \mid \mathcal{F}_t \ ].
\end{eqnarray}
By Jensen's inequality, this leads to
\begin{eqnarray*}
    &&\hspace{-1cm} \| E[\vartheta_{t+1:w}(q) - \vartheta_{t+1:w}(h)\mid \mathcal{F}_t \ ] \|_p  \\
        &\leq& \| q_{t+1} - h_{t+1} \|_p + \| E[\vartheta_{t+2:w-1}(q)(1 - \theta_{f,t+1}(h)) - \vartheta_{t+2:w-1}(h) \theta_{f,t+1}^-(h) \mid \mathcal{F}_t \ ] \|_p \\
        &=& \| q_{t+1} - h_{t+1} \|_p +
            \| E[(\vartheta_{t+2:w-1}(q) - \vartheta_{t+2:w-1}(h)) \theta_{f,t+1}^-(h) \mid \mathcal{F}_t \ ] \|_p \\
        &\leq& \| q_{t+1} - h_{t+1} \|_p + \| E[\vartheta_{t+2:w-1}(q) - \vartheta_{t+2:w-1}(h) \mid \mathcal{F}_{t+1} \ ] \|_p.
\end{eqnarray*}
The proof is completed by induction.
\end{proof}

\begin{proof}[Proof of Proposition \ref{prop:first err decomp}]
Introduce the regression function
\begin{equation}\label{equ:reg function rec}
    \bar{\rho}_{\mathcal{H},t}(x) = E[\vartheta_{t+1:w}(f,\hat{q}_{\mathcal{H}}) \mid X_t = x]
\end{equation}
of $\vartheta_{t+1:w}(f,\hat{q}_{\mathcal{H}})$ and let
\begin{equation}\label{equ:target regress}
    \bar{q}_{\mathcal{H}, t} = \pi_{\mathcal{H}_t} \bar{\rho}_{\mathcal{H},t} =
        \pr_{\mathcal{H}_t} \vartheta_{t+1:w}(f,\hat{q}_{\mathcal{H}})
\end{equation}
be its projection onto $\mathcal{H}_t$. By the triangle inequality
\begin{eqnarray}\label{equ:first err decomp 1}
    \|\hat{q}_{\mathcal{H},t} - q_t \|_2 &\leq&
        \|\hat{q}_{\mathcal{H},t} - \bar{q}_{\mathcal{H},t} \|_2 +
        \|\bar{q}_{\mathcal{H},t} - \bar{\rho}_{\mathcal{H},t} \|_2 +
        \|\bar{\rho}_{\mathcal{H},t} - q_t \|_2.
\end{eqnarray}
Again by the triangle inequality and because $\mathcal{H}_t$ is convex so that the projection
$\pi_{\mathcal{H}_t}$ from $L_2(\mathbb{R}^m,\mu_t)$ onto
$\mathcal{H}_t$ is distance decreasing
\begin{eqnarray}\label{equ:first err decomp 1a}
    \|\bar{q}_{\mathcal{H},t} - \bar{\rho}_{\mathcal{H},t} \|_2
    &=& \|\pi_{\mathcal{H}_t} \bar{\rho}_{\mathcal{H},t} - \bar{\rho}_{\mathcal{H},t} \|_2  \nonumber \\
    &\leq& \|\pi_{\mathcal{H}_t} \bar{\rho}_{\mathcal{H},t} - \pi_{\mathcal{H}_t} q_t \|_2 +
        \|\pi_{\mathcal{H}_t} q_t - q_t \|_2 + \| q_t - \bar{\rho}_{\mathcal{H},t} \|_2  \\
    &\leq& \|\pi_{\mathcal{H}_t} q_t - q_t \|_2 + 2 \| q_t - \bar{\rho}_{\mathcal{H},t} \|_2. \nonumber
\end{eqnarray}
Inserting \eqref{equ:first err decomp 1a} back into \eqref{equ:first err decomp 1} gives
\begin{equation}\label{equ:first err decomp 1b}
    \|\hat{q}_{\mathcal{H},t} - q_t \|_2 \leq \inf_{h \in \mathcal{H}_t} \|h - q_t \|_2 +
        \|\hat{q}_{\mathcal{H},t} - \bar{q}_{\mathcal{H},t} \|_2 +
        3 \|\bar{\rho}_{\mathcal{H},t} - q_t \|_2.
\end{equation}
By Lemma \ref{lem:err dec convex}
\begin{equation}\label{equ:first err decomp 1c}
    \|\hat{q}_{\mathcal{H},t} - \bar{q}_{\mathcal{H},t} \|_2 =
        \|\hat{q}_{\mathcal{H},t} - \pi_{\mathcal{H}_t} \bar{\rho}_{\mathcal{H},t} \|_2 \leq
        E[l_t(\hat{q}_{\mathcal{H}})]^{1/2}.
\end{equation}
For the third term in \eqref{equ:first err decomp 1b}, by Proposition
\ref{prop:continuity z functional}
\begin{equation}\label{equ:first err decomp LS}
    \|\bar{\rho}_{\mathcal{H},t} - q_t \|_2 =
        \| E[ \vartheta_{t+1:w}(f,\hat{q}_{\mathcal{H}}) - \vartheta_{t+1:w}(f,q) \mid X_t \ ] \|_2  \leq
        \sum_{s=t+1}^{t+w+1} \|\hat{q}_{\mathcal{H},s} - q_s \|_2.
\end{equation}
\end{proof}

\subsection{Covering Number Bounds}
\label{sec:covering number bounds}

We define the so called centered loss class
\begin{equation}\label{equ:centered loss class}
    \mathcal{L}_t(\mathcal{H}) = \{l_t(h) \mid h \in \mathcal{H} \}.
\end{equation}
To bound the fluctuations of the sample error
$E[l_t(\hat{q}_{\mathcal{H}})]^{1/2}$ later on in Section
\ref{sec:proof consistency thm} we require bounds on the empirical
$L_1$-covering numbers $N(\eps, \mathcal{L}_t(\mathcal{H}),
d_{1,P_n})$ of the centered loss class.

The first step is to bound the covering numbers of
$\mathcal{L}_t(\mathcal{H})$ in terms of the covering numbers of
$\mathcal{H}_t$ and the cash flow class which is defined as
\begin{equation}\label{equ:class G}
    \mathcal{G}_t = \{\vartheta_{t+1:w}(f,h) \mid h \in \mathcal{H}\}.
\end{equation}
\begin{lem}\label{lem:centered loss class entropy bound}
Let $1 \leq p \leq \infty$. If $\mathcal{H}_t $ is uniformly
bounded by $H$ and the cash flow class $\mathcal{G}_t$ by $\Theta$
then for $w \geq 0$
\begin{equation}\label{equ:centered loss class cov bound 1}
    N(8(H + \Theta) \eps, \mathcal{L}_t(\mathcal{H}), d_{p,P_n}) \leq
        N\left(\eps,\mathcal{H}_t, d_{p,P_n}\right)^2 N\left(\eps,\mathcal{G}_t,
        d_{p,P_n}\right)^2.
\end{equation}
For $w=0$ the estimate \eqref{equ:centered loss class cov bound 1}
simplifies to
\begin{equation}\label{equ:centered loss class cov bound 1 w=0}
    N(8(H + \Theta) \eps, \mathcal{L}_t(\mathcal{H}), d_{p,P_n}) \leq
        N\left(\eps,\mathcal{H}_t, d_{p,P_n}\right)^2 N\left(\eps,\mathcal{H}_{t+1}, d_{p,P_n}\right)^2.
\end{equation}
\end{lem}

Note that if the payoff functions $f$ is in
$L_{\infty}(\mathbf{X})$ and the approximation spaces
$\mathcal{H}_t$ are uniformly bounded by $H$ then
$\vartheta_{t+1:w}(f,h) \leq \Theta \equiv \max(\|f\|_{\infty},
H)$ and the assumptions of Lemma \ref{lem:centered loss class
entropy bound} are satisfied.

\begin{proof}
We first recall some basic properties of covering numbers. If
$\mathcal{F}$ and $\mathcal{G}$ are two classes of functions and
$\mathcal{F} \pm \mathcal{G} = \{f \pm g \mid f \in \mathcal{F}, g
\in \mathcal{G}\}$ is the class of formal sums or differences,
then for all $1 \leq p \leq \infty$
\begin{equation}\label{equ:cov num sum}
    N(\eps, \mathcal{F} \pm \mathcal{G}, d_{p,P_n}) \leq N\left(\frac{\eps}{2}, \mathcal{F}, d_{p,P_n}\right)
        N\left(\frac{\eps}{2}, \mathcal{G}, d_{p,P_n}\right).
\end{equation}
Furthermore, if $\mathcal{G}$ class of functions uniformly bounded
by $G$, it follows from $\|g_1^2 - g_2^2\|_{p,P_n}^p = P_n
(g_1-g_2)^p(g_1+g_2)^p \leq (2G)^p\|g_1 - g_2\|_{p,P_n}^p$ that
\begin{equation}\label{equ:cov num square}
    N(\eps, \mathcal{G}^2, d_{p,P_n}) \leq N\left(\frac{\eps}{2 G}, \mathcal{G}, d_{p,P_n}\right),
\end{equation}
Enlarging a class increases the covering numbers. Now
\begin{equation}\label{equ:centered loss class expr}
    \mathcal{L}_t(\mathcal{H}) \subset (\mathcal{H}_t - \mathcal{G}_t)^2 - (\pr_{\mathcal{H}_t}\mathcal{G}_t - \mathcal{G}_t)^2.
\end{equation}
Because $\pr_{\mathcal{H}_t} \mathcal{G}_t \subset \mathcal{H}_t$,
it is sufficient to bound the covering number of the slightly
larger class
\begin{equation}\label{equ:centered loss class red expr}
    \tilde{\mathcal{L}}_t(\mathcal{H}) = (\mathcal{H}_t - \mathcal{G}_t)^2 - (\mathcal{H}_t - \mathcal{G}_t)^2.
\end{equation}
If $\mathcal{H}_t$ is uniformly bounded by $H < \infty$ and
$\vartheta_{t+1,w}(f,h) \leq \Theta$, we get from \eqref{equ:cov
num sum} and \eqref{equ:cov num square}
\begin{equation}\label{equ:centered loss class cov bound 1bis}
    N(\eps, \tilde{\mathcal{L}}_t(\mathcal{H}), d_{p,P_n}) \leq
        N\left(\frac{\eps}{8(H + \Theta)},\mathcal{H}_t, d_{p,P_n}\right)^2 N\left(\frac{\eps}{8(H + \Theta)},\mathcal{G}_t, d_{p,P_n}\right)^2.
\end{equation}
For $w=0$ the Lipschitz bound \eqref{equ:continuity TvR} directly
leads to
\begin{equation}\label{equ:cov bound target van Roy}
    N(\eps,\mathcal{G}_t, d_{p,P_n}) \leq N(\eps,\mathcal{H}_{t+1},
    d_{p,P_n}).
\end{equation}
\eqref{equ:centered loss class cov bound 1 w=0} follows directly
from \eqref{equ:centered loss class cov bound 1bis} and
\eqref{equ:cov bound target van Roy}.
\end{proof}

A simple example for which tight covering number bounds exists are
subsets of linear vector spaces. If $\mathcal{H}_t = \{h \in
\mathcal{K} \mid \|h\|_{\infty} \leq R\}$ and $\mathcal{K}$ is a
linear vector space of dimension $d$ then
\begin{equation}\label{equ:cov number linear space subset}
    N(\eps, \mathcal{H}_t, d_{2,P_n}) \leq N(\eps, \{h \in \mathcal{K} \mid P_n h^2 \leq R^2\}, d_{2,P_n})
        \leq \left(\frac{4R + \eps}{\eps}\right)^d.
\end{equation}
The first inequality in \eqref{equ:cov number linear space subset}
is obvious because $\mathcal{H}_t$ is a subset of $\{h \in
\mathcal{K} \mid P_n h^2 \leq R^2\}$. The second inequality is
standard and can be found for instance in
\citeasnoun{Carl-Stephani:1990} or
\citeasnoun{Vaart-Wellner:1996}.

\eqref{equ:cov number linear space subset} would provide uniform
covering number estimates for \eqref{equ:centered loss class cov
bound 1 w=0} in case of linear approximation spaces and $w = 0$.
We can not apply \eqref{equ:cov number linear space subset} to
upper bound the right hand side of \eqref{equ:centered loss class
cov bound 1} in the general situation $w > 0$ because the cash
flow class $\mathcal{G}_t$ is not anymore a subset of a linear
space, even if the underlying approximation space $\mathcal{H}_t$
is a finite dimensional linear vector space. This is where the
Vapnik-Chervonenkis theory comes into play.

An important type of function classes for which good uniform
estimates on the covering numbers exist without assuming any
linear structure are the so called Vapnik-Chervonenkis classes or
VC-classes, introduced in \citeasnoun{Vapnik-Chervonenkis:1971}
for classes of indicator functions, i.e., classes of sets. Let
$\mathcal{C}$ be a class of subsets of a set $S$. We say that the
class $\mathcal{C}$ picks out a subset $A$ of a set $\sigma_n =
\{x_1,\ldots, x_n\} \subset S$ of $n$ elements if $A = C \cap
\sigma_n$ for some $C \in \mathcal{C}$. The class $\mathcal{C}$ is
said to shatter $\sigma_n$ if each of its $2^n$ subset can be
picked out by $\mathcal{C}$. The VC-dimension of $\mathcal{C}$ is
the largest integer $n$ such that there exists a set of $n$ points
which can be shattered by $\mathcal{C}$, i.e.,
\begin{equation}\label{equ:VC dim}
    \vc(\mathcal{C}) = \sup \{ n \mid \Delta_n(\mathcal{C}) = 2^n\},
\end{equation}
where
\begin{equation}\label{equ:VC dim a}
    \Delta_n(\mathcal{C}) = \max_{\{x_1,\ldots, x_n\}} \card \{ C \cap \{x_1,\ldots, x_n\} \mid C \in \mathcal{C} \}
\end{equation}
is the so called growth or shattering function. A class
$\mathcal{C}$ is called a Vapnik-Chervonenkis or VC-class if
$\vc(\mathcal{C}) < \infty$. A VC-class of dimension $d$ shatters
no set of $d+1$ points. The ``richer'' the class $\mathcal{C}$ is,
the larger the cardinality of sets which still can be shattered.
We illustrate it by a simple example. The class of left open
intervals $\{ (-\infty, c] \mid c \in \mathbb{R}\}$ cannot shatter
any two-point set because it cannot pick out the largest of the
two points and therefore has VC-dimension one. By similar
reasoning, the class of intervals $\{ (-a, b] \mid a,b \in
\mathbb{R}\}$ shatters two-point sets but fails to shatter
three-point sets: it cannot pick out the largest and the smallest
point of a three-point set. Contrary, the collection of closed
convex subsets of $\mathbb{R}^2$ has infinite VC-dimension:
Consider a set $\sigma_n$ of $n$ points on the unit circle. Every
subset $A \subset \sigma_n$ of the $2^n$ subsets can be picked out
by the closed convex hull $\overline{\text{co}}(A)$ of $A$. A
peculiar property of a VC-class is that the shattering function of
VC-classes grows only polynomially in $n$, more precisely we have
the following result which is due to Sauer, Vapnik-Chervonenkis
and Shelah, see \citeasnoun{Vaart-Wellner:1996}, Corollary 2.6.3,
or \citeasnoun{Dudley:1999}.
\begin{lem}[Sauer's Lemma]\label{lem:Sauer}
If $\mathcal{C}$ is a VC-class with VC-dimension $d =
\vc(\mathcal{C})$, then
\begin{equation}\label{equ:Sauer}
    \Delta_n(\mathcal{C}) \leq \sum_{i=0}^d {n \choose i} \leq 1.5 \frac{n^d}{d!} \leq \Big(\frac{e\;
    n}{d}\Big)^d.
\end{equation}
\end{lem}
VC-classes have a variety of permanence properties which allow the
construction of new VC-classes from basic VC-classes by simple
operations such as complements, intersections, unions or products.
We again refer to \citeasnoun[section 2.6.5]{Vaart-Wellner:1996},
or \citeasnoun{Dudley:1999}.

The concept of VC-classes of sets can be extended to classes of
functions in several ways. A common approach is to associate to a
class of functions its subgraph class. More precisely, the
subgraph of a real-valued function $g$ on an arbitrary set $S$ is
defined as
\begin{equation}\label{equ:subgraph}
    \mathrm{Gr}(g) = \{(x,t) \in S \times \mathbb{R} \mid t \leq g(x)\}.
\end{equation}
A class of real-valued functions $\mathcal{G}$ on $S$ is called a
VC-subgraph class, or just VC-class, if its class of subgraphs is
a VC-class and the VC-dimension of $\mathcal{G}$ is defined as
\begin{equation}\label{equ:subgraph vc dim}
    \vc(\mathcal{G}) = \vc(\{\mathrm{Gr}(g) \mid g \in \mathcal{G} \}).
\end{equation}
An equivalent definition is obtained by extending the notion of
shattering. A class of real-valued functions $\mathcal{G}$ is said
to shatter a set $\{x_1,\ldots, x_n\} \subset S$ if there is $r
\in \mathbb{R}^n$ such that for every $b \in \{0,1\}^n$, there is
a function $g \in \mathcal{G}$ such that for each $i$, $g(x_i) >
r_i$ if $b_i = 1$, and $g(x_i) \leq r_i$ if $b_i = 0$. The
definition
\begin{equation}\label{equ:subgraph vc dim bis}
    \vc(\mathcal{G}) = \sup \{n \mid \exists \{x_1,\ldots, x_n\} \subset S \text{ shattered by } \mathcal{G}\}
\end{equation}
agrees with \eqref{equ:subgraph vc dim}. For the proof note that a
set is shattered by the subgraph class $\{\mathrm{Gr}(g) \mid g
\in \mathcal{G} \}$ if and only if it is shattered by the class of
indicator functions $\{\theta(g(x) - t) \mid g \in \mathcal{G}\}$,
where $\theta(s) = 1_{\{s \geq 0\}}$. The VC-dimension
\eqref{equ:subgraph vc dim bis} for classes of functions is often
called pseudo-dimension, see \citeasnoun{Pollard:1990} and
\citeasnoun{Haussler:1995}. An alternative generalization is
obtained by so called VC-major classes, originally introduced by
Vapnik. For more details on the relation of the two concepts we
refer to \citeasnoun{Dudley:1999}.

\begin{lem}\label{lem:VC dim vector space}
Let $\mathcal{G}$ be a finite dimensional real vector space of
measurable real-valued functions. Then, the class of sets
$\mathcal{G}^+ = \{\{g \geq 0\} \mid g \in \mathcal{G}\}$ is a VC
class with $\vc(\mathcal{G}^+) \leq \dim(\mathcal{G})$. If $g_0$
is a fixed function, then  $\vc((g_0 + \mathcal{G})^+) =
\vc(\mathcal{G}^+)$. Finally, $\mathcal{G}$ is a VC-class and
$\vc(\mathcal{G}) = \dim(\mathcal{G})$.
\end{lem}

\begin{proof}
For the first two statements we refer to \citeasnoun[theorem
4.2.1]{Dudley:1999}, or \citeasnoun[section
2.6]{Vaart-Wellner:1996}. The last statements follows from the
first two: Let $g_0(x,t) = -t$ and consider the affine class of
functions $g_0 + \mathcal{G}$ on $S \times \mathbb{R}$. Then, the
subgraph class of $\mathcal{G}$ is precisely $(g_0 +
\mathcal{G})^+$.
\end{proof}

An important property of VC-classes is that their covering numbers
$N(\eps, \mathcal{G}, d_{p,\mu})$ are polynomial in $\eps^{-1}$
for $\eps \rightarrow 0$. More precisely we have the following
estimates for the covering numbers of VC-classes due to
\citeasnoun{Haussler:1995}, see also \citeasnoun[Theorem
2.6.7]{Vaart-Wellner:1996}.

\begin{lem}\label{lem:Haussler cov num bound}
Let $\mathcal{G} \subset L_p(\mu)$ be a class of functions with an
envelope $G \in L_p(\mu)$, i.e., $g \leq G$ for all $g \in
\mathcal{G}$. Then,
\begin{equation}\label{equ:Haussler cov num bound}
    N(\eps \|G\|_{p,\mu}, \mathcal{G}, d_{p,\mu}) \leq e (\vc(\mathcal{G})+1) 2^{\vc(\mathcal{G})} \left(\frac{2 e
    }{\eps}\right)^{p \,\vc(\mathcal{G})}.
\end{equation}
\end{lem}

After this short digression on VC-theory we continue estimating
the empirical $L_1$-covering numbers of the centered loss class
$\mathcal{L}_t(\mathcal{H})$. The next result is fundamental to
generalize the estimate \eqref{equ:centered loss class cov bound 1
w=0} to a strictly positive look-ahead parameter $w > 0$. It
bounds the VC-dimension of $\mathcal{G}_t$ in terms of the
VC-dimension of the approximation spaces
$\mathcal{H}_{t+1},\ldots, \mathcal{H}_{t+w+1}$.

\begin{prop}\label{prop:VC dim bound target LS}
Assume that for all $s \geq t$, $\mathcal{H}_s$ are VC-classes of
functions with $\vc(\mathcal{H}_s) \leq d$. Then $\mathcal{G}_t$
is a VC-class with VC-dimension
\begin{equation}\label{equ:VC dim bound target LS}
    \vc(\mathcal{G}_t) \leq c(w) d,
\end{equation}
where $c(w) = 2(w+2)\log_2(e(w+2))$.
\end{prop}

Inequalities \eqref{equ:centered loss class cov bound 1},
\eqref{equ:centered loss class cov bound 1 w=0},
\eqref{equ:Haussler cov num bound}, and \eqref{equ:VC dim bound
target LS} finally lead to explicit uniform bounds for the
empirical $L_1$-covering numbers of the centered loss class
$\mathcal{L}_t(\mathcal{H})$.

\begin{cor}\label{cor:covering number bounds}
Assume that all $\mathcal{H}_s$ are classes of function uniformly
bounded by $H$ and with bounded VC-dimension $\vc(\mathcal{H}_s)
\leq d$. If the cash flow function satisfies
$\vartheta_{t+1:w}(f,h) \leq H$, then
\begin{eqnarray}\label{equ:covering num bound}
    N(\eps, \mathcal{L}_t(\mathcal{H}), d_{1,P_n}) &\leq&  \\[2mm]
        && \hspace{-2.5cm} \left \{\begin{array}{ll}
          e^4(d+1)^2(c(w)d+1)^2 \left(\displaystyle{\frac{64 e H}{\eps}}\right)^{2d(c(w)+1)}, & \text{ for } \; w \geq 1, \\[3mm]
          e^4(d+1)^4 \left(\displaystyle{\frac{64 e H}{\eps}}\right)^{4d},                    & \text{ for } \; w = 0.\\[2mm]
        \end{array} \right. \nonumber
\end{eqnarray}
\end{cor}

Optimal stopping is a particular stochastic control problem with a
simple control space. The proof of Proposition \ref{prop:VC dim
bound target LS} relies on the observation that the VC-dimension
of the class of indicator functions $\mathcal{C}_s =
\{\theta_{f,s}(h) \mid h_s \in \mathcal{H}_s\}$, which appear in
the definition of $\tau_t(h)$ and $\vartheta_{t+1:w}(h)$, is
bounded by $\vc(\mathcal{H}_s)$. It is an interesting question how
Proposition \ref{prop:VC dim bound target LS} can be extended to
more general stochastic control problems.

Before we proceed to the proof of Proposition \eqref{prop:VC dim
bound target LS} we add a remark on VC-classes and their
VC-dimension. Let $\mathcal{A}$ be a class of sets. The class of
indicator functions $\{1_A \mid A \in \mathcal{A}\}$ is a VC-class
in the sense that its subgraph class is a VC-class if and only if
$\mathcal{A}$ is a VC-class and $\vc(\mathcal{A}) = \vc(\{1_A \mid
A \in \mathcal{A}\})$. Let $\theta(x) = 1_{\{x \geq 0\}}$. If
$\mathcal{A}$ is a VC-class, $\vc(\mathcal{A})=d$, then by Sauer's
Lemma \ref{lem:Sauer}, for $x_1,\ldots, x_n$ and all $t \in
\mathbb{R}^n$
\begin{equation}\label{equ:VC class vs VC subgraph class}
    \card \{(\theta(1_A(x_i) - t_i))_{i=1,\ldots,n} \mid A \in
    \mathcal{A}\} \leq \Big(\frac{e\; n}{d}\Big)^{d}.
\end{equation}
Conversely, if we find a polynomial bound like \eqref{equ:VC class
vs VC subgraph class}, $\mathcal{A}$ must be a VC-class and we can
bounds its VC-dimension.

To prove Proposition \ref{prop:VC dim bound target LS} we first
establish the following general result on VC-classes.

\begin{lem}\label{lem:VC 2 step}
Let $\mathcal{X}$, $\mathcal{Y}$ be two sets and $\mathcal{A}$,
$\mathcal{B}$ VC-classes of subsets of $\mathcal{X}$
(respectively, $\mathcal{Y}$). Assume that $\vc(\mathcal{A}) \leq
d$, $\vc(\mathcal{B}) \leq d$. Let $f : \mathcal{X} \rightarrow
\mathbb{R}$ and $g : \mathcal{Y} \rightarrow \mathbb{R}$ be
non-negative functions. Define the class of functions
\begin{equation}\label{equ:VC 2 step def}
    \mathcal{F}(\mathcal{A}, \mathcal{B}) = \{F_{A,B}(x,y) =
    1_A(x) f(x) + 1_{A^c}(x)1_B(y)g(y) \mid A \in \mathcal{A}, B \in \mathcal{B} \}
\end{equation}
Then $\mathcal{F}(\mathcal{A}, \mathcal{B})$ is a VC-subgraph
class, its growth function is bounded by
\begin{equation}\label{equ:VC 2 step growth estim}
    \Delta_n(\mathcal{F}(\mathcal{A}, \mathcal{B})^+) \leq
        \left(\frac{e\;n}{d}\right)^{2d},
\end{equation}
and
\begin{equation}\label{equ:VC 2 step vc dim}
    \vc(\mathcal{F}(\mathcal{A}, \mathcal{B})) \leq 2 d \log_2(e).
\end{equation}
The estimates \eqref{equ:VC 2 step growth estim} and \eqref{equ:VC
2 step vc dim} generalize to
\begin{equation}\label{equ:VC 2 step def H}
    \mathcal{F}(\mathcal{A}, \mathcal{H}) = \{F_{A,h}(x,y) =
    1_A(x) f(x) + 1_{A^c}(x)h(y) \mid A \in \mathcal{A}, h \in \mathcal{H} \}
\end{equation}
where $\mathcal{H}$ is a VC-class of function with
$\vc(\mathcal{H}) = \vc(\mathcal{H}^+) \leq d$.
\end{lem}

\begin{proof}
Given points $(x_i, y_i) \in \mathcal{X} \times \mathcal{Y}$ and
$t_i \in \mathbb{R}$, $i = 1, \ldots, n$, we need to bound the
cardinality of
\begin{equation}\label{equ:VC 2 step 1}
    \{(\theta(F_{A,B}(x_i,y_i) - t_i))_{i=1,\ldots,n} \mid A \in \mathcal{A}, B \in
    \mathcal{B}\},
\end{equation}
as a subset of the binary cube $\{0,1\}^n$. Because
\begin{equation*}
    F_{A,B}(x,y) = 1_B(y_i)(g(y_i) - 1_A(x_i)g(y_i)) +
    1_A(x_i)f(x_i),
\end{equation*}
and $(g(y_i) - 1_A(x_i)g(y_i)) \geq 0$ we find that
\begin{multline}\label{equ:VC 2 step 2}
    \theta(F_{A,B}(x_i,y_i) - t_i) = \\ \left \{
        \begin{array}{ll}
          \theta(1_B(y_i) - \tau_i(A)) & \text{on } S_+(A) = \{(x_j,y_j) \mid 1_{A^c}(x_j)g(y_j) > 0 \} \\[1mm]
          \theta(1_A(x_i)f(x_i) - t_i) & \text{on } S_0(A) = \{(x_j,y_j) \mid 1_{A^c}(x_j)g(y_j) = 0 \} \\
        \end{array}
        \right . ,
\end{multline}
where
\begin{equation}\label{equ:VC 2 step 3}
    \tau_i(A) = \frac{t_i - 1_A(x_i) f(x_i)}{g(y_i) - 1_A(x_i)g(y_i)}.
\end{equation}
Fix $A$ and vary $B$ over $\mathcal{B}$. Because $\vc(\mathcal{B})
\leq d$ we see from \eqref{equ:VC 2 step 2} and Sauer's lemma,
that the binary set
\begin{equation}\label{equ:VC 2 step 4}
   \{(\theta(F_{A,B}(x_i,y_i) - t_i))_{i=1,\ldots,n} \mid B \in
    \mathcal{B}\}
\end{equation}
has cardinality $K$ bounded above by $\left(e n d^{-1}\right)^d$.
Let
\begin{equation}\label{equ:VC 2 step 5}
   b_1(A), \ldots, b_K(A)
\end{equation}
enumerate the distinct elements of \eqref{equ:VC 2 step 4}
generated by sets $B_k$. For $(x_i,y_i) \in S_0(A)$ we have
\begin{equation}\label{equ:VC 2 step 6}
    b_{k,i}(A) = \theta(1_A(x_i)f(x_i) - t_i),
\end{equation}
and if $(x_i,y_i) \in S_+(A)$
\begin{multline}\label{equ:VC 2 step 7}
    b_{k,i}(A) = \theta(1_{B_k}(y_i) - \tau_i(A)) = \\ \left \{
        \begin{array}{ll}
          \theta(1_A(x_i) - \tau_i(B_k)) & \text{on } S_+(B_k) = \{(x_j,y_j) \mid f(x_j) - 1_{B}(x_j)g(y_j) > 0 \} \\[1mm]
          1 - \theta(1_A(x_i) - \tau_i(B_k)) &  \text{on } S_-(B_k) = \{(x_j,y_j) \mid f(x_j) - 1_{B}(x_j)g(y_j) < 0 \}\\[1mm]
          \theta(1_{B_k}(y_i)g(y_i) - t_i) & \text{on } S_0(B_k) = \{(x_j,y_j) \mid f(x_j) - 1_{B}(x_j)g(y_j) = 0 \}\\
        \end{array}
        \right . .
\end{multline}
Consequently, Sauer's lemma again implies that for each fixed $k$
the binary set
\begin{equation}\label{equ:VC 2 step 8}
    \{b_{k}(A) \mid A \in \mathcal{A} \}
\end{equation}
has cardinality at most $\left(e n d^{-1}\right)^d$. This proves
\eqref{equ:VC 2 step growth estim}. Again by Sauer's lemma, very
$n_0 > 0$ such that
\begin{equation}\label{equ:Pollard dim bound card}
    \card\{(\theta(F_{A,B}(x_i,y_i) - t_i))_{i=1,\ldots,n} \mid A \in \mathcal{A}, B \in
    \mathcal{B}\} \leq \left( \frac{e n}{d}\right)^{2d} < 2^n,
\end{equation}
for all $n > n_0$ is an upper bound of
$\vc(\mathcal{F}(\mathcal{A}, \mathcal{B})^+)$. To find $n_0$, we
look for solutions $n_0 = d j$ that are multiples of $d$.
\eqref{equ:Pollard dim bound card} leads to the condition
\begin{equation*}
    \log_2(e j) < j,
\end{equation*}
which is satisfied for example by $j = 2  \log_2(e)$. The
extension to $\mathcal{F}(\mathcal{A}, \mathcal{H})$ is
straightforward. Replace $\theta(1_B(y_i) - \tau_i(A))$ in
\eqref{equ:VC 2 step 2} by $\theta(h(y_i) - \tau_i(A))$, where
$\tau_i(A) = (t_i - f(x_i))/1_{A^c}(x_i)$ and follow the same
lines of reasoning.
\end{proof}

\begin{proof}[Proposition \ref{prop:VC dim bound target LS}]
Recall definition \eqref{equ:z ipol} of the cash flow function,
according to which
\begin{equation}\label{equ:def target bis}
    \vartheta_{t+1:w}(f,h) = \theta_{f,t+1}(h)f_{t+1} + \ldots +
        \theta_{f,t+w+1}(h) \prod_{r=t+1}^{t+w}\theta_{f,r}^-(h) f_{t+w+1}
        + \prod_{r=t+1}^{t+w+1}\theta_{f,r}^-(h) h_{t+w+1}.
\end{equation}
Because the classes of indicator functions
\begin{equation}\label{equ:theta class}
    \mathcal{C}_s = \{\theta_{f,s}(h) = 1_{\{f_s-h_s \geq 0\}} \mid h_s \in \mathcal{H}_s\}, \quad
    \mathcal{C}^-_s = \{\theta^-_{f,s}(h) = 1_{\{f_s-h_s < 0\}} \mid h_s \in \mathcal{H}_s\},
\end{equation}
are VC classes with VC-dimension
\begin{equation}\label{equ:theta class 1}
    \vc(\mathcal{C}^-_s) = \vc(\mathcal{C}_s) = \vc((f_s - \mathcal{H}_s)^+) =
    \vc(\mathcal{H}_s^+) = \vc(\mathcal{H}_s) \leq d,
\end{equation}
we can recursively apply Lemma \ref{lem:VC 2 step} to derive the
bound
\begin{equation}\label{equ:proof VC dim bound target LS}
    \card \{(\theta(\vartheta_{t+1:w}(f,h)(\mathbf{x}_i) - t_i))_{i=1,\ldots,n} \mid h \in \mathcal{H}\}
        \leq \left( \frac{e n}{d}\right)^{d (w+2)}.
\end{equation}
The VC-dimension of $\mathcal{G}_t$ is then estimated as in the
proof of Lemma \ref{lem:VC 2 step}. This completes the proof of
Proposition \ref{prop:VC dim bound target LS}.
\end{proof}

\subsection{Proof of Theorem \ref{thm:ERM consistency} and Theorem \ref{thm:ERM conv rate}}
\label{sec:proof consistency thm}

The centered loss $l_t(\hat{q}_{\mathcal{H}})$ depends on the
sample $D_n$. To control the fluctuations of the random variable
$E[l_t(\hat{q}_{\mathcal{H}})]$ we need uniform estimates over the
whole centered loss class $\mathcal{L}_t(\mathcal{H})$. The usual
procedure is to apply exponential deviation inequalities for the
empirical process
\begin{equation}\label{equ:emp process}
    \{\sqrt{n}(E[l] - P_n l) \mid l \in \mathcal{L}_t(\mathcal{H}) \}
\end{equation}
indexed by $\mathcal{L}_t(\mathcal{H})$, which are closely related
to Uniform Law of Large Numbers. For background we refer to
\citeasnoun{Pollard:1984}, \citeasnoun{Vaart-Wellner:1996},
\citeasnoun{Talagrand:1994}, and \citeasnoun{Kohler-etal:2002}.

The application of standard deviation inequalities to the whole
centered loss class $\mathcal{L}_t(\mathcal{H})$ is not efficient
since the empirical minimizer is close to the actual
$L_2$-minimizer with high probability. Therefore, the random
element $l_t(\hat{q}_{\mathcal{H}})$ is with high probability in a
small subset of $\mathcal{L}_t(\mathcal{H})$. To get sharper
estimates, the empirical process needs to be localized such that
more weight is assigned to these loss functions.
\citeasnoun{Lee-Bartlett:1996} proved the following localized
deviation inequality.

\begin{thm}[\citeasnoun{Lee-Bartlett:1996}, Theorem 6]\label{thm:Lee etal}
Let $\mathcal{L}$ be a class of functions such that $|l| \leq
K_1$, $E[l] \geq 0$, and for some $K_2 \geq 1$,
\begin{equation}\label{equ:cond L}
    E[l^2] \leq K_2 E[l] \quad \forall l \in \mathcal{L}.
\end{equation}
Let $a, b > 0$ and $0 < \delta < \half$. Then, for all
\begin{equation}\label{equ:cond on n}
    n \geq \max\left( \frac{4(K_1 + K_2)}{\delta^2(a + b)}, \frac{K_1^2}{\delta^2(a + b)}\right),
\end{equation}
\begin{eqnarray}\label{equ:Lee etal}
    \mathbb{P}\left(
        \sup_{l \in \mathcal{L}}
                \frac{E[l] - P_n(l)}{E[l] + a + b} \geq \delta
    \right) &\leq& \nonumber \\
    &&\hspace{-3cm} 2 \sup_{x_1,\ldots,x_{2n} \in \mathcal{X}^{2n}} N\left(\frac{\delta b}{4}, \mathcal{L}, d_{1,P_{2n}}\right)
            \exp\left(-\frac{3 \delta^2 a n }{4K_1 + 162 K_2}\right)  + \\
    &&\hspace{-3cm} 4 \sup_{x_1,\ldots,x_{2n} \in \mathcal{X}^{2n}} N\left(\frac{\delta b}{4K_1}, \mathcal{L}, d_{1,P_{2n}}\right)
            \exp\left(-\frac{\delta^2 a n }{2K_1^2}\right), \nonumber
\end{eqnarray}
where $P_{2n}$ is the empirical measure supported at
$(x_1,\ldots,x_{2n})$.
\end{thm}
A similar bound has been obtained by \citeasnoun[Proposition
7]{Cucker-Smale:2001} for $L_{\infty}$-covering numbers. Theorem
\ref{thm:Lee etal} has been improved in \citeasnoun{Kohler:2000}
by applying chaining techniques, and in
\citeasnoun{Bousquet-Bartlett-Mendelson:2002a} by using
concentration properties of local Rademacher averages. For
additional background on related bounds we refer to
\citeasnoun{Talagrand:1994}, \citeasnoun{Ledoux:1996},
\citeasnoun{Massart:2000}, and \citeasnoun{Rio:2001}. The
advantage of Theorem \ref{thm:Lee etal}, as compared to the
Pollard's deviation inequality, is that it improves the quadratic
dependence on $\eps$ in standard deviation inequalities to a
linear dependence.

The centered loss has a special structure which allows to bound
its variance in terms of its expectation.
\begin{lem}\label{lem:Bernstein type}
Let $\mathcal{H}_t$ be convex, uniformly bounded by $H < \infty$,
and assume that $\vartheta_{t+1:w}(f,h) \leq \Theta$ for some
constant $\Theta < \infty$. Then the centered loss class
$\mathcal{L}_t(\mathcal{H})$ is uniformly bounded and for all $l
\in \mathcal{L}_t(\mathcal{H})$
\begin{eqnarray}\label{equ:Bernstein type}
    |l|    &\leq& 4H(\Theta + H), \nonumber \\
    E[l^2] &\leq& 4(\Theta + H)^2 E[l].
\end{eqnarray}
\end{lem}

\begin{proof} We get from the definition \eqref{equ:centered loss}
of $l_t(h)$ that
\begin{eqnarray}\label{equ:Bernstein type 2}
    l_t(h) &=& (h_t - \pr_{\mathcal{H}_t} \vartheta_{t+1:w}(f,h))
        (h_t + \pr_{\mathcal{H}_t} \vartheta_{t+1:w}(f,h) - 2 \vartheta_{t+1:w}(f,h))  \\
     &\leq& 2(\Theta + H) (h_t - \pr_{\mathcal{H}_t} \vartheta_{t+1:w}(f,h)).\nonumber
\end{eqnarray}
Therefore,
\begin{equation*}
    E[l_t(h)^2] \leq 4(\Theta + H)^2 E[|h_t - \pr_{\mathcal{H}_t} \vartheta_t|^2]
        \leq 4(\Theta + H)^2 E[l_t(h)],
\end{equation*}
where the last step follows form Lemma \ref{lem:err dec convex}.
\end{proof}

Our plan is to apply Theorem \ref{thm:Lee etal} to a suitably
scaled loss class
\begin{equation}\label{equ:scaled loss class}
    \lambda \mathcal{L}_t(\mathcal{H}) = \{ \lambda l \mid l \in
    \mathcal{L}_t(\mathcal{H})\},
\end{equation}
where we choose $\lambda$ such that $|\lambda l| \leq 1$ (the
scaling gives a term $\beta_n^2$ in the consistency condition
\eqref{equ:ERM consistency conv in prob cond} instead of
$\beta_n^4$). Because an empirical risk minimizer satisfies
$P_n(l_t(\hat{q}_{\mathcal{H}})) \leq 0$, it follows that for any
$\eps > 0$ and scaling factor $\lambda > 0$
\begin{eqnarray}\label{equ:err decomp conv rate 3}
    \mathbb{P}\left( E[l_t(\hat{q}_{\mathcal{H}})] \geq \eps \right) &\leq&
         \mathbb{P}\left( E[l_t(\hat{q}_{\mathcal{H}})] \geq 2P_n(l_t(\hat{q}_{\mathcal{H}})) + \eps \right)  \nonumber \\
         \nothing{ &=& \mathbb{P}\left(
            \frac{E[l_t(\hat{q}_{\mathcal{H}})] - P_n(l_t(\hat{q}_{\mathcal{H}}))}{E[l_t(\hat{q}_{\mathcal{H}})] + \eps}
            \geq \frac{1}{2} \right) \\}
         &\leq& \mathbb{P}\left(
            \sup_{l \in \lambda \mathcal{L}_t(\mathcal{H})}
                \frac{E[l] - P_n(l)}{E[l] + \lambda \eps} \geq \frac{1}{2}
                \right).
\end{eqnarray}
Assume that the conditions of Lemma  \ref{lem:Bernstein type} are
satisfied and set $\beta = \max(\Theta, H)$. If we choose the
scaling factor $\lambda = 1/(8\beta^2)$ the scaled class $\lambda
\mathcal{L}_t(\mathcal{H})$ satisfies
\begin{eqnarray}\label{equ:Bernstein type rescaled}
    |\lambda l|    &\leq& 1, \nonumber \\
    E[(\lambda l)^2] &\leq& 2 E[\lambda l].
\end{eqnarray}
Theorem \eqref{thm:Lee etal} applied with $\mathcal{L} = \lambda
\mathcal{L}_t(\mathcal{H})$, $K_1 = 1$, $K_2=2$, $a = b =
\eps/(16\beta^2)$, $\delta = 1/2$ implies
\begin{equation}\label{equ:err decomp conv rate 4}
    \mathbb{P}\left( E[l_t(\hat{q}_{\mathcal{H}})] \geq \eps \right)
         \leq 6\sup_{x_1,\ldots,x_{2n}}
            N\left(\frac{\eps}{128\beta^2}, \frac{1}{8\beta^2} \mathcal{L}_t(\mathcal{H}), d_{1,P_{2n}}\right)
            \exp\left(-\frac{n \eps}{6998 \beta^2}\right) ,
\end{equation}
for $n \geq 382\beta^2/\eps$. The $\eps$-covering number of
$\lambda \mathcal{L}_t(\mathcal{H})$ is the same as the
$(\lambda^{-1} \eps)$-covering number of the unscaled class
$\mathcal{L}_t(\mathcal{H})$. If the VC-dimension of
$\mathcal{H}_s$, $s \geq t$ are bounded by $d$ the covering number
bound \eqref{equ:covering num bound} shows that
\begin{eqnarray}\label{equ:main probability tail estimate}
    \mathbb{P}\left( E[l_t(\hat{q}_{\mathcal{H}})] \geq \eps \right)
         &\leq& K \left(\frac{1}{\eps}\right)^v \exp\left(-\frac{n
         \eps}{6998 \beta^2}\right),
\end{eqnarray}
where
\begin{equation}\label{equ:err decomp conv rate 5}
    v = v(w,d) = 2d(c(w) + 1),
\end{equation}
and
\begin{equation}\label{equ:err decomp conv rate 6}
    K = K(d, w, \beta) = 6 e^4(d+1)^2(c(w)d+1)^2(1024 e \beta)^{v(d,w)}.
\end{equation}

\begin{proof}[Proof of Theorem \ref{thm:ERM consistency}]
$\beta_n$ is a sequence of truncation thresholds tending to
infinity. If $\bar{q}_{\beta_n,t}$ is the continuation value for
the truncated payoff $T_{\beta_n} f$ we get from
\eqref{equ:truncation conv} that $\|q_t - \bar{q}_{\beta_n,t}\|_2
\rightarrow 0$. The error decomposition \eqref{equ:first err
decomp} separates the approximation error and the sample error.
The denseness assumption implies that the approximation error
$\inf_{h \in \mathcal{H}_{n,t}} \|h - \bar{q}_{\beta_n,t} \|_2$
tends to zero if $n \rightarrow \infty$. It remains to analyze the
sample error $E[l_t(\hat{q}_{\mathcal{H}_n})]$ for underlying
payoff $T_{\beta_n} f$.

We apply \eqref{equ:main probability tail estimate} to
$\mathcal{H} = \mathcal{H}_n$ for which $d = d_n$ and $\beta =
\beta_n$. There exists a constant $C(\eps, w)$ such that for every
fixed $\eps > 0$
\begin{equation}\label{equ:proof ERM constistency 1}
    \mathbb{P} (E[l_t(\hat{q}_{\mathcal{H}_n})] \geq \eps)
        \leq C(\eps,w) \exp \left(d_n \log(\beta_n) - \frac{n \eps}{6998
        \beta_n^2}\right).
\end{equation}
The right hand side converges to zero for every fixed $\eps > 0$
if $n/\beta_n^2$ diverges to infinity faster than $d_n
\log(\beta_n)$ or if $d_n \beta_n^2 \log(\beta_n) n^{-1}
\rightarrow 0$. Convergence in probability follows from
\eqref{equ:first err decomp} by induction. Convergence in
$L_1(\mathbb{P})$ is shown by evaluating
\begin{equation}\label{equ:proof ERM constistency 3}
    \mathbb{E}[E[l_t(\hat{q}_{\mathcal{H}_n})]] \leq \eps + \int_{\eps}^{\infty} \mathbb{P}(E[l_t(\hat{q}_{\mathcal{H}_n})] > t) dt,
\end{equation}
using the estimate \eqref{equ:proof ERM constistency 1}.
Conditions \eqref{equ:ERM consistency conv in prob cond} and
\eqref{equ:ERM consistency conv almost sure cond} imply
\begin{eqnarray}\label{equ:proof ERM constistency 2}
    \sum_{n=1}^{\infty} \mathbb{P} (E[l_t(\hat{q}_{\mathcal{H}_n})] \geq \eps)
        &\leq& C(\eps,w) \sum_{n=1}^{\infty} \exp \left(d_n \log(\beta_n) - \frac{n \eps}{6998 \beta_n^2}\right) \nonumber \\
        &=& C(\eps,w) \sum_{n=1}^{\infty} n^{-\frac{n}{\log(n)\beta_n^2 }
            \left(\frac{\eps}{6998} - \frac{d_n \beta_n^2 \log(\beta_n)}{n}
            \right)} < \infty.
\end{eqnarray}
Almost sure convergence follows from the Borel-Cantelli Lemma.
\end{proof}

\begin{proof}[Proof of Theorem \ref{thm:ERM conv rate}]
Integrating \eqref{equ:main probability tail estimate} over $\eps$
shows that for any $\kappa \geq \frac{1}{n}$
\begin{eqnarray}\label{equ:err decomp conv rate concent 6}
    \mathbb{E}\left[E[l_t(\hat{q}_{\mathcal{H}})]\right] &=&
         \int_0^{\infty} \mathbb{P}\left( E[l_t(\hat{q}_{\mathcal{H}})] \geq \eps \right) d\eps \nonumber \\
        &\leq& \kappa + K n^v \int_{\kappa}^{\infty} \exp\left(-\frac{n \eps}{6998 \beta^2}\right) d\eps \nonumber \\
        &\leq& \kappa + K n^{v-1} 6998 \beta^2 \exp\left(-\frac{n \kappa}{6998 \beta^2}\right).
\end{eqnarray}
Setting
\begin{equation}\label{equ:err decomp conv rate 7}
    \kappa = \frac{6998 \beta^2}{n}\log\left(6998 K \beta^2 n^v\right),
\end{equation}
leads to the upper bound
\begin{equation}\label{equ:err decomp conv rate 8}
    \mathbb{E}\left[E[l_t(\hat{q}_{\mathcal{H}})]\right] \leq
        \frac{6998 \beta^2 + \log(6998 K \beta^2)}{n} + \frac{v
        \log(n)}{n}.
\end{equation}
Corollary \ref{cor:first err decomp} implies that
%
%
\begin{equation*}
    \mathbb{E}\left[\|\hat{q}_{\mathcal{H},t} - q_t \|_2^2\right]
        \leq 2 \cdot 16^{w+1} \left( \max_{s=t,\ldots,t+w+1} \; \inf_{h \in \mathcal{H}_{s}} \|h - q_s \|_2^2 +
              \; \mathbb{E}\left[\max_{s=t,\ldots,t+w+1}E[l_s(\hat{q}_{\mathcal{H}})]\right] \right).
\end{equation*}
But
\begin{equation*}
    \mathbb{E}\left[\max_{s=t,\ldots,t+w+1}E[l_s(\hat{q}_{\mathcal{H}})]\right] \leq
        (w+2) \max_{s=t,\ldots,t+w+1} \mathbb{E}\left[E[l_s(\hat{q}_{\mathcal{H}})]\right].
\end{equation*}
Apply \eqref{equ:err decomp conv rate 8} to complete the proof.
\end{proof}

\begin{proof}[Proof of Corollary \ref{cor:sample complexity bound}]
Estimate \eqref{equ:main probability tail estimate} implies
\begin{equation}\label{equ:sample complexity 1}
    \mathbb{P}\left( E[l_t(\hat{q}_{\mathcal{H}})] \geq \epsilon \right) \leq
        K \exp\left( -\frac{n\epsilon}{13996 \beta^2}\right) \exp\left( -\frac{n\epsilon}{13996 \beta^2} - \log(\epsilon) v \right).
\end{equation}
By straightforward calculations, the right hand side is smaller
than $\delta$ for all $n$ satisfying
\begin{equation}\label{equ:sample complexity 2}
    n \geq 13996 \beta^2 \max\left(\frac{1}{\epsilon} \log\left( \frac{K}{\delta}\right), v \log\left( \frac{1}{\epsilon}\right) \right).
\end{equation}
The sample complexity bound \eqref{equ:sample complexity bound 1}
follows from Corollary \ref{cor:first err decomp} and
\eqref{equ:sample complexity 1}, \eqref{equ:sample complexity 2}
with $\epsilon = \eps / (32 (w+2) 16^{w})$.

\end{proof}

\subsection{Proof of Corollary \ref{cor:Sobolev approx rates}}
\label{sec: proof Sobolev approx rates}

Because $q_t \in W^k(L_{\infty}(I,\lambda))$, Jackson type
estimates imply that for every $r > k$ there exists a polynomial
$p_r \in \mathcal{P}_r$
\begin{equation}\label{equ:Sobolev approx rates}
    \|p_r - q_t \|_{\infty,I,\lambda} \leq
        C_I \, r^{-k} \, \| q_t \|_{\infty,k,I,\lambda} .
\end{equation}
The constant $C_I$ only depends on $I$ but not on $r$ or $q_t$.
See for instance \citeasnoun[Theorem 6.2, Chapter
7]{DeVore-Lorentz:1993}. Consequently
\begin{equation}\label{equ:Sobolev approx rates 1}
    \|p_r \|_{\infty,I,\lambda} \leq \|p_r - q_t
    \|_{\infty,I,\lambda} + \| q_t \|_{\infty,I,\lambda}
        \leq 2 \, \| q_t \|_{\infty,k,I,\lambda}
\end{equation}
for $r$ sufficiently large. We therefore may restrict the
minimization to the convex, uniformly bounded set of functions
$\mathcal{H}_{n,t}$ as defined in \eqref{equ:Sobolev approx rates
approx space}. The VC-dimension of $\mathcal{H}_{n,t}$ is bounded
by $n^{m/(m+2k)}$. Theorem \ref{thm:ERM conv rate} applies.
Because $X_t$ is localized to $I$, the approximation error in
\eqref{equ:ERM conv rate error estimate} is bounded by
\begin{equation}\label{equ:Sobolev approx rates 2}
    \inf_{p \in \mathcal{H}_{n,t}} \|p - q_t \|_{2}^2 \leq
        \inf_{p \in \mathcal{H}_{n,t}} \|p - q_t \|_{\infty,I,\lambda}^2 \leq
        C_I \, n^{-2k/(2k+m)} \, \| q_t \|_{\infty,k,I,\lambda} .
\end{equation}
Inserting $\vc(\mathcal{H}_{n,t}) \leq n^{m/(m+2k)}$ into
\eqref{equ:ERM conv rate error estimate} shows that the sample
error is of the order $$O\left(\log(n) n^{-2k/(2k + m)}\right).$$
The extension to $\mu_t$ with bounded density with respect to
Lebesgue measure is proved identically. \hfill \qed

\subsection{Proof of Proposition \ref{prop:truncation}}
\label{sec: proof trunc error estimate}

Note that
\begin{equation}\label{equ:max Lipschitz}
    |\max(a,x) - \max(a,y)| \leq |x-y|.
\end{equation}
The representation of the continuation value in terms of the
transition functions gives
\begin{eqnarray*}
    \| q_t - \bar{q}_{\beta,t} \|_p
        &=& \| E[\max(f_{t+1},q_{t+1})) - \max(T_{\beta} f_{t+1},\bar{q}_{\beta,t+1}) \mid X_t ] \|_p  \\
        &\leq& \| E[f_{t+1} - T_{\beta} f_{t+1} \mid X_t ] \|_p +
            \| E[q_{t+1} - \bar{q}_{\beta,t+1} \mid X_t ]  \|_p  \\
        &\leq& \| (f_{t+1}-\beta) 1_{\{f_{t+1} > \beta\}} \|_p + \| q_{t+1} - \bar{q}_{\beta,t+1} \|_p.
\end{eqnarray*}
If $f \in L_p(\mathbf{X})$, then $\| (f_{t+1}-\beta) 1_{\{f_{t+1}
> \beta\}} \|_p \rightarrow 0$ for $\beta \rightarrow \infty$. We first
recall that for a nonnegative random variable $Y$ and $r > 1$
\begin{equation}\label{equ:Lp norm equation}
    E[Y^r] = r \int_{0}^{\infty} y^{r-1} P(Y > y)  \, d y.
\end{equation}
Then \eqref{equ:truncation rate} follow from
\begin{eqnarray*}
    \|(f_{t+1}-\beta) 1_{\{f_{t+1} > \beta\}}\|_r^{r} &=&
        r \int_0^{\infty} u^{r-1} P((f_{t+1}-\beta) 1_{\{f_{t+1} > \beta\}} > u) \, d u \\
        &=& r \int_{\beta}^{\infty} (u-\beta)^{r-1} P(f_{t+1} > u) \, d u \\
        &\leq& r \int_{\beta}^{\infty} u^{r-1} P(f_{t+1}^p > u^p) \, d u \\
        &\leq&  \frac{r}{p-r} E[f_{t+1}^p] \beta^{r-p} \leq O(\beta^{r-p}),
\end{eqnarray*}
where we have used Markov's inequality to get to the last line.
\hfill \qed

\nocite{Dudley:1999}
\nocite{Vapnik:1999}
\nocite{Pollard:1984}
\nocite{Vaart-Wellner:1996}

%
\bibliographystyle{kluwer}

\end{document}